%% file: conj-LQYZ-r.tex
\documentclass{amsart}

\usepackage{delimset, amssymb, enumitem, mathtools}
\setlist{itemsep=4pt, topsep=0pt, leftmargin=17pt, listparindent=11pt}

\usepackage{xcolor}
\definecolor{BIT}{cmyk}{1, 0, 1, 0}
\usepackage[colorlinks, citecolor=BIT, ocgcolorlinks]{hyperref}
\usepackage{aliascnt}

\usepackage[capitalize]{cleveref}
\crefname{conjecture}{Conjecture}{Conjectures}
\crefname{itm}{}{}
\Crefname{itm}{}{}
\creflabelformat{itm}{~\upshape#2#1#3}
\usepackage{subcaption}
\AddToHook{env/conjecture/begin}{\crefalias{theorem}{conjecture}}
\AddToHook{env/lemma/begin}{\crefalias{theorem}{lemma}}
\AddToHook{env/proposition/begin}{\crefalias{theorem}{proposition}}
\usepackage{booktabs}
\newtheorem{theorem}{Theorem}[section]

\newtheorem{lemma}[theorem]{Lemma}
\newtheorem{proposition}[theorem]{Proposition}
\numberwithin{equation}{section}
\numberwithin{table}{section}

\allowdisplaybreaks

\usepackage[numbers, sort&compress, nonamebreak, merge, elide, longnamesfirst]{natbib}
\usepackage{etoolbox}
\usepackage{hyphenat}
\makeatletter
\patchcmd{\NAT@citexnum}
{\@citea \NAT@test{\@ne}\NAT@spacechar\NAT@mbox{\NAT@super@kern\NAT@@open}}
{\@citea \nohyphens{\NAT@test{\@ne}}\nobreakspace{}\NAT@mbox{\NAT@super@kern\NAT@@open}}
{}
{\PackageError{nonbreaking-citet}
{Failed to patch \string\citet}
{The internal definition of natbib may have changed.}}
\makeatother

\usepackage{tikz}
\tikzset{
  edge/.style={semithick},
  ball/.style={shape=circle, minimum size=1mm, ball color=black, inner sep=0.5},
  border/.style={very thick}, 
  ellipsis/.style={shape=circle, fill, inner sep=.5}, 
  ribbon/.style={black, thick}, 
  round/.style={line cap=round, line join=round}, 
  vertex/.style={circle, fill=black, inner sep=1.5pt}}

\DeclareMathOperator{\hgt}{ht}
\DeclareMathOperator{\inc}{inc}
\DeclareMathOperator{\sgn}{sgn}
\DeclareMathOperator{\stc}{sc}
\DeclareMathOperator{\sh}{sh}

\newlength{\casewidth}
\title{Exact thresholds for Schur positivity of the lattices $\mathbf m\times\mathbf 2$ and $\mathbf m\times\mathbf 3$}

\author[D.G.L.~Wang]{David G.L. Wang$^*$}
\address[David G.L. Wang]{School of Mathematics and Statistics \& MIIT Key Laboratory of Mathematical Theory and Computation in Information Security, Beijing Institute of Technology, Beijing 102400, P.\ R.\ China.}
\email{glw@bit.edu.cn}
\thanks{$^*$Wang is the corresponding author, and is supported by the NSFC (Grant No.\ 12171034).}

\author[K. Zhang]{K. Zhang}
\address[K. Zhang]{School of Mathematics and Statistics, Beijing Institute of Technology, Beijing 102400, P.\ R.\ China.}
\email{kai@bit.edu.cn}



\begin{document}

\begin{abstract}
We determine the exact thresholds for Schur positivity in the two families of chain products $\mathbf m\times\mathbf2$ and $\mathbf m\times\mathbf3$: the former is Schur positive exactly for $m\le7$, and the latter exactly for $m\le6$. For $m\ge8$, we prove non-Schur-positivity in both families by exhibiting explicit negative Schur coefficients obtained from Pieri's rules and stable-composition counts. The remaining finite cases are settled by exact SageMath computations; in particular, $\mathbf7\times\mathbf3$ has a negative Schur coefficient. These results settle the $n=2$ and $n=3$ cases in the conjectural picture of Li, Qiu, Yang, and Zhang and sharpen the $n=3$ boundary by one. We also show that $\mathbf m\times\mathbf3$ is not strongly nice for $m\ge44$.
\end{abstract}

\maketitle

\section{Introduction}

George \citet{Bir1912} introduced the \emph{chromatic polynomial} of a graph $G$ in his study of the Four Color Problem more than a century ago. It is now a fundamental graph invariant and remains a central object of study in algebraic combinatorics. \citet{Sta95} generalized it to the \emph{chromatic symmetric function}:
\[
X_G=\sum_{\kappa\colon V(G)\to\{1, 2, \dots\}}x_1^{\abs{\kappa^{-1}(1)}}x_2^{\abs{\kappa^{-1}(2)}}\dotsm, 
\]
where the sum runs over proper colorings $\kappa$ of $G$. 

A major line of research on chromatic symmetric functions concerns positivity. Given a basis $\{b_\lambda\}$ of the vector space $\Lambda$ of symmetric functions, a symmetric function is \emph{$b$-positive} if all its coefficients in this basis are nonnegative. The Schur basis is one of the most prominent bases of $\Lambda$; see \citet{Mac95B}. Schur positivity connects combinatorics, representation theory, and geometry because Schur coefficients often record multiplicities of irreducible representations; see \citet{Sta81,Sta99B,Sag01B}. The basis of elementary symmetric functions $\{e_\lambda\}$ is also important, and every $e$-positive symmetric function is Schur positive \cite{Sta99B}. 

A graph $G$ is \emph{$b$-positive} if $X_G$ is $b$-positive. \citet{Sta95} reformulated a conjecture of \citet{SS93} as the assertion that the chromatic symmetric function of the incomparability graph of every $(3+1)$-free poset is $e$-positive. The weaker conclusion of Schur positivity was established by \citet{Gas96}. The stronger assertion, known as the \emph{Stanley--Stembridge conjecture}, was recently proved by \citet{Hik25}. Typical examples include unit interval graphs, also called claw-free interval graphs, such as complete graphs, paths, lollipops, and $K$-chains; see \citet{Ale21,Dv18,GS01,SW16}. \citet[p.~187]{Sta95} also asked, ``Which $X_G$ is $e$-positive?'' Considerable progress has been made on this question; see \citet{WZ25,TV26-vglue} and the references therein.

Another major conjecture, posed by \citet{Sta98} and attributed there to \citeauthor{Gas99}, asserts that every claw-free graph is Schur positive; see also \cite{Gas99}. This conjecture remains open. Among results not obtained through $e$-positivity, \citet{Sv25} proved Schur positivity for generalized net graphs. For graph families that may contain claws, \citet{WW20} characterized the Schur-positive complete bipartite and tripartite graphs, and \citet{Sv26X} extended this classification to all complete multipartite graphs. \citet{TW23X} established Schur positivity for the spiders $S(a,2,1)$ and $S(a,4,1)$.

A complementary line of work seeks explicit obstructions to Schur positivity. Building on the special ribbon tabloids of \citet{ER90}, \citet{WW20} gave a combinatorial formula for every Schur coefficient of $X_G$ and used it to produce negative coefficients; see \cref{thm:scoeff}. A partition~$\lambda$ \emph{dominates} a partition $\mu$ if $\sum_{i=1}^k\lambda_i\ge \sum_{i=1}^k\mu_i$ for every $k$; we write $\lambda\ge\mu$. Following \citet{Sta98}, a graph $G$ is \emph{nice} if, whenever $\lambda\ge\mu$, the existence of a stable partition of type $\lambda$ implies the existence of one of type $\mu$. Equivalently,
\[
[m_\lambda]X_G>0\Longrightarrow[m_\mu]X_G>0, 
\]
where $m_\lambda$ denotes the monomial symmetric function. \citet{Sta98} established the following implication.

\begin{proposition}[\citeauthor{Sta98}]\label{prop:Sta98}
If $G$ is Schur positive, then it is nice.
\end{proposition}

\citet{LLYZ25} defined a graph $G$ to be \emph{strongly nice} if
\[
[m_\lambda]X_G\le[m_\mu]X_G
\]
whenever $\lambda\ge\mu$ in dominance order. They proved that every Schur-positive graph is strongly nice and used this stronger condition to show that no squid graph is Schur positive, thereby confirming a conjecture of \citet{WW20}.

These notions extend naturally to posets: a poset $P$ is Schur positive (resp., strongly nice or nice) if its incomparability graph $\inc(P)$ has the corresponding property. Distributive lattices form a central class of finite posets. By Birkhoff's representation theorem, every finite distributive lattice is isomorphic to the lattice of order ideals of its poset of join-irreducible elements; under this representation, joins and meets correspond to unions and intersections; see \citet{BirGarrett1940B}.

The study of Schur positivity for distributive lattices originates in a longstanding problem about Boolean algebras. A conjecture of \citet{Gri98} is equivalent to the assertion that the incomparability graph $\inc(B_n)$ of the Boolean algebra $B_n$ is nice. Motivated by this conjecture and \cref{prop:Sta98}, \citet[p.~270]{Sta98} observed that ``$\inc(B_n)$ might be $s$-positive, which is true for $n\leq 4$.'' Because Boolean algebras form a special class of distributive lattices, he further speculated, ``Perhaps even the incomparability graph of any distributive lattice is $s$-positive.'' These observations motivated subsequent work on Schur positivity for distributive lattices. 

Recently, \citet{LQYZ25} answered Stanley's broad question in the negative by constructing distributive lattices that are not nice and hence not Schur positive. A basic family of distributive lattices consists of products $\mathbf m\times\mathbf n$ of two chains. These products are distributive and, by \citet{LE99}, always nice. Nevertheless, \citet{LQYZ25} proved that the following subfamily is not Schur positive.

\begin{theorem}[\citeauthor{LQYZ25}]\label{thm:LQYZ}
For any integers $k\ge 5$ and $n\ge (k+2)/2$, the distributive lattice $\mathbf{(n+k)}\times\mathbf n$ is not Schur positive.
\end{theorem}

\citet{LQYZ25} proved this theorem by showing that the $s_\rho$-coefficient in the Schur expansion of~$X_{\inc(\mathbf{(n+k)}\times\mathbf n)}$ is negative, where 
\[
\rho=(2n+k-1,\, 2n+k-3,\, \dots,\, k+3,\, k-3,\, 2,\, 2)\vdash n(n+k)
\]
has length $n+2$. 
As an extension of \cref{thm:LQYZ},
\citet[Conjecture~5.1]{LQYZ25} proposed a conjectural
non-Schur-positive range for products of two chains. The printed formulation
of the conjecture is ambiguous: it gives the two conditions
\[
(m\ge 8 \text{ and } n\ge 2)
\qquad\text{or}\qquad
(m\ge n+5 \text{ and } n\ge 3).
\]
Read literally, the second condition is contained in the first. Moreover,
the first condition includes the equal-factor cases $m=n\ge8$, whereas
\citet[Conjecture~5.2]{LQYZ25} predicts Schur positivity for products of
equal chains. For the two cases considered in this manuscript, however, there is
no ambiguity: Conjecture~5.1 predicts that $\mathbf m\times\mathbf2$ and
$\mathbf m\times\mathbf3$ are not Schur positive for every $m\ge8$. The
SageMath computations reported immediately after Conjecture~5.1 begin at
$m=8$ for $n=2$, at~$m=8=n+5$ for $n=3$, and at $m=9=n+5$ for $n=4$.

Writing $m=n+k$, the conditions $k\ge5$ and
$n\ge(k+2)/2$ in \cref{thm:LQYZ} are equivalent to
\[
n+5\le m\le3n-2.
\]
In particular, \cref{thm:LQYZ} applies only when $n\ge4$.
Since the lattices $\mathbf m\times\mathbf n$ and
$\mathbf n\times\mathbf m$ are isomorphic, we may assume
$m\ge n$ without loss of generality. In particular, the cases
$n=2$ and $n=3$ are not covered by \cref{thm:LQYZ}.
For $n\ge4$, after imposing $m\ge n$, this theorem leaves
the two complementary ranges
\[
n\le m\le n+4
\qquad\text{and}\qquad
m\ge3n-1.
\]
The Schur-positivity status in these ranges remains open in
general, although finite computational evidence is available
for some values.
\citet{LQYZ25} reported SageMath computations
showing non-Schur-positivity for
\[
8\le m\le1199\quad\text{when }n=2,
\qquad
8\le m\le100\quad\text{when }n=3,
\]
and
\[
9\le m\le50\quad\text{when }n=4, 
\]
which supplied substantial finite evidence.

In this manuscript, we determine the exact thresholds for the two smallest
nontrivial values $n=2,3$:
\[
\mathbf m\times\mathbf2
\text{ is not Schur positive if and only if }m\ge8,
\]
whereas
\[
\mathbf m\times\mathbf3
\text{ is not Schur positive if and only if }m\ge7.
\]
Consequently, we prove the range predicted by \citeauthor{LQYZ25}
for $n=2,3$; for $n=3$, our result improves the conjectured threshold
from $m\ge8$ to $m\ge7$.

The values $n=2,3$ are the two smallest nontrivial cases and are precisely those not covered by \cref{thm:LQYZ}. They are also structurally distinct: the incomparability graph is bipartite for $n=2$, whereas $n=3$ is the first non-bipartite case. The partitions and arguments of \citet{LQYZ25} are formulated for $n\ge4$ and do not directly extend to these two cases.

The main challenge is to select partitions that index negative Schur coefficients of $X_{\inc(\mathbf m\times\mathbf n)}$. Computation can verify individual values of $m$, but it does not supply a uniform target partition or a closed formula. The two isolated vertices in the incomparability graph introduce an additional constraint: after writing
\[
X_{\inc(\mathbf m\times\mathbf n)}=s_1^2X_{H_m^n},
\]
we must determine which coefficients of $X_{H_m^n}$ can contribute. Pieri's rules organize this selection. We use target partitions of length $3$ for $n=2$ and length $4$ for $n=3$, which keeps the ensuing stable-composition enumeration tractable.

For the target partitions used below, Pieri's rules and the
independence-number bound leave only a short signed list of special
ribbon tabloids. The corresponding Schur coefficient is therefore a
difference of stable-composition counts. In both the $n=2$ and
$n=3$ cases, the contributions associated with negative-sign
tabloids eventually dominate those associated with positive-sign
tabloids. The resulting formulas are
\[
-\frac{m(m-2)(m-7)}{3}
\]
for $\mathbf m\times\mathbf2$ and
\[
-\frac43
\Bigl((m-8)^3+12(m-8)^2+44(m-8)+24\Bigr)
\]
for the uniform $\mathbf m\times\mathbf3$ calculation. Thus the
obstruction detected by our argument is a signed imbalance among the
stable-composition counts that survive the Pieri and
independence-number filters.

We also obtain two byproducts. First, the graphs $G_m$, which are closely related to the half-graphs introduced by Erd{\H o}s and Hajnal, are not Schur positive for every $m\ge6$; see \cref{prop:s-coeff:half}. Second, \cref{thm:sn-:m3} shows that $\mathbf m\times\mathbf3$ is not strongly nice for $m\ge44$. Thus we obtain an infinite family of distributive lattices that fail even the stronger necessary condition for Schur positivity.

\section{Preliminaries}\label{sec:pre}

Let $n$ be a positive integer and write $[n]=\{1,\dots,n\}$. A \emph{composition} of $n$ is a sequence of positive integers~$\lambda_1,\dots,\lambda_\ell$ with sum $n$, written~\mbox{$\lambda=\lambda_1\dotsm\lambda_\ell\vDash n$}. Its \emph{length} is $\ell(\lambda)=\ell$. A \emph{partition} of $n$ is a composition with weakly decreasing parts, written $\lambda=\lambda_1\dotsm\lambda_\ell\vdash n$. Exponents indicate repeated parts; for example, $(4,4,2,2,1)=(4^2,2^2,1)$. We identify a partition $\lambda$ with its Young diagram $Y_\lambda$, whose $i$th row has $\lambda_i$ left-justified cells. We use the French convention, with rows indexed from bottom to top. \Cref{ydop} depicts the Young diagram of $(4^2,2^2,1)$.

\begin{figure}[h]
\centering
\begin{tikzpicture}[scale=0.4]
\draw[very thick] (0,0) -- (4,0)--(4,2)--(2,2)--(2,4)--(1,4)--(1,5)--(0,5)--(0,0);
\foreach \x in {1, 2, 3} {
\draw[gray] (\x, 0) -- (\x, 2);}
\draw[gray] (0,1)--(4,1);
\draw[gray] (1,2)--(1,5);
\draw[gray] (0,2)--(2,2);
\draw[gray] (0,3)--(2,3);
\draw[gray] (0,4)--(2,4);
\end{tikzpicture}
\caption{The Young diagram of shape $(4^2,2^2,1)$.}
\label{ydop}
\end{figure}	

If $\lambda_i\ge\mu_i$ for all $i$, the \emph{skew shape} $\lambda/\mu$ is $Y_\lambda\setminus Y_\mu$. A \emph{horizontal strip} (resp., \emph{vertical strip}) is a skew shape with no two cells in the same column (resp., row). A \emph{ribbon} is a connected skew shape containing no $2\times2$ square. For a ribbon~$R$, its \emph{size}~$\abs{R}$ is its number of cells, its \emph{height}~$\hgt(R)$ is one less than the number of rows, its \emph{sign} is
\[
\sgn(R)=(-1)^{\hgt(R)}, 
\]
and its \emph{head} is the leftmost cell on the top row. A \emph{ribbon tabloid} of shape $Y$ is a decomposition of $Y$ into ribbons. The empty tabloid is special. Recursively, a nonempty ribbon tabloid is a \emph{special ribbon tabloid} if every ribbon head lies in the first column and deleting the top ribbon leaves a smaller special ribbon tabloid. \Cref{ribbon} gives examples of shape $(4^2,2^2,1)$. \Citet{ER90} introduced these objects under the name \emph{special rim hook tabloids}. For a special ribbon tabloid $T$, its \emph{shape} $\sh(T)$ is $Y$, its \emph{content} $\kappa_T$ is the composition of ribbon sizes read from top to bottom, and its \emph{type} $\lambda_T$ is the partition obtained by rearranging the parts of $\kappa_T$. Once $\sh(T)$ is fixed, $T$ is determined by $\kappa_T$. Its \emph{height} is $\hgt(T)=\sum_{R\in T}\hgt(R)$, and its \emph{sign} is the product of the ribbon signs. For example, the special ribbon tabloid in the right-hand panel of \cref{ribbon} has height $2$, content $(3, 2, 5, 3)$, and sign $1$.

\begin{figure}[h]
\centering
\begin{subfigure}[b]{0.4\textwidth}
\centering
\captionsetup{labelformat=empty}
\begin{tikzpicture}[scale=0.4]
\draw[very thick] (0,0) -- (4,0)--(4,2)--(2,2)--(2,4)--(1,4)--(1,5)--(0,5)--(0,0);
\foreach \x in {1, 2, 3} {
\draw[gray] (\x, 0) -- (\x, 2);}
\draw[gray] (0,1)--(4,1);
\draw[gray] (1,2)--(1,5);
\draw[gray] (0,2)--(2,2);
\draw[gray] (0,3)--(2,3);
\draw[gray] (0,4)--(2,4);
\foreach \x/\y in {0.5/0.5, 3.5/0.5, 3.5/1.5, 0.5/2.5, 0.5/3.5, 0.5/4.5, 0.5/1.5, 1.5/3.5}{
\fill (\x, \y) circle (3pt);}
\draw[thick] (1.5, 3.5) -- (1.5, 1.5) -- (3.5, 1.5);
\draw[thick] (0.5, 4.5) -- (0.5, 3.5);
\draw[thick] (0.5, 2.5) -- (0.5, 1.5);
\draw[thick] (0.5, 0.5) -- (3.5, 0.5);
\end{tikzpicture}
\caption{}
\end{subfigure}
\begin{subfigure}[b]{0.4\textwidth}
\centering
\captionsetup{labelformat=empty}
\begin{tikzpicture}[scale=0.4]
\draw[very thick] (0,0) -- (4,0)--(4,2)--(2,2)--(2,4)--(1,4)--(1,5)--(0,5)--(0,0);
\foreach \x in {1, 2, 3} {
\draw[gray] (\x, 0) -- (\x, 2);}
\draw[gray] (0,1)--(4,1);
\draw[gray] (1,2)--(1,5);
\draw[gray] (0,2)--(2,2);
\draw[gray] (0,3)--(2,3);
\draw[gray] (0,4)--(2,4);
\foreach \x/\y in {0.5/4.5, 1.5/3.5, 0.5/2.5, 1.5/2.5, 0.5/1.5, 3.5/0.5, 0.5/0.5, 2.5/0.5}{
\fill (\x, \y) circle (3pt);}
\draw[thick] (0.5, 4.5) -- (0.5, 3.5) -- (1.5, 3.5);
\draw[thick] (0.5, 2.5) -- (1.5, 2.5);
\draw[thick] (0.5, 1.5) -- (3.5, 1.5) -- (3.5, 0.5);
\draw[thick] (0.5, 0.5) -- (2.5, 0.5);
\end{tikzpicture}
\caption{}
\end{subfigure}
\caption{Left: a ribbon tabloid. Right: a special ribbon tabloid.}
\label{ribbon}
\end{figure}

Let $G=(V,E)$ be a graph of order $n$. Following \citet{WW20}, a \emph{semi-ordered partition} of $G$ is a set partition of $V$ in which the blocks of each fixed cardinality are linearly ordered. A \emph{composition} of $G$ is a sequence $\pi=(B_1,\dots,B_\ell)$ of pairwise disjoint nonempty subsets whose union is~$V$. It has \emph{type} $\tau=(\tau_1,\dots,\tau_\ell)\vDash n$ if $\abs{B_i}=\tau_i$ for all $i$. For fixed $\tau$, semi-ordered partitions with block sizes given by the parts of $\tau$ are naturally equivalent to type-$\tau$ compositions: for each size $r$, place the ordered blocks of size $r$ in the positions $i$ for which $\tau_i=r$. We use the language of compositions throughout.

A subset $B\subseteq V(G)$ is \emph{stable} if the induced subgraph $G[B]$ has no edges. The \emph{independence number}~$\alpha(G)$ is the order of a largest stable set of $G$. A composition $\pi$ of $G$ is said to be \emph{stable} if each block of $\pi$ is stable. We denote by $\mathcal{SC}(G;\tau)$ the set of type-$\tau$ stable compositions, namely,
\[
\begin{aligned}
\mathcal{SC}(G;\tau)
=\{(B_1,\dots,B_{\ell(\tau)})\colon{}&
(B_1,\dots,B_{\ell(\tau)})\text{ is a composition of $G$},\\
& B_i\text{ is stable and }\abs{B_i}=\tau_i\text{ for all $i$}\}.
\end{aligned}
\]
Its cardinality is written 
\[
\stc(G;\tau)=\abs{\mathcal{SC}(G;\tau)}.
\]
$\stc(G;\tau)$ depends only on the partition $\lambda_\tau$ obtained by rearranging the parts of $\tau$; thus 
\[
\stc(G;\tau)=\stc(G;\lambda_\tau).
\]  
When the graph $G$ and the composition $\tau$ are clear from context, we suppress them from the notation and write simply $\mathcal{SC}$ and $\stc$.

We write $[s_\mu]X_G$ for the coefficient of $s_\mu$ in $X_G$. The following formula is due to Wang and Wang; see~\cite{WW20}.

\begin{theorem}[\citeauthor{WW20}]\label{thm:scoeff}
For any graph $G$ and any integer partition $\mu$,
\[
[s_{\mu}]X_G
=\sum_{T\in \mathcal{T}_{\mu}}
\sgn(T)
\stc(G; \lambda_T), 
\]
where $\mathcal{T}_{\mu}$ is the set of special ribbon tabloids of shape $\mu$. 
\end{theorem}

For a basis $\{b_\lambda\}_{\lambda\vdash n}$ of the degree-$n$ symmetric functions, a graph $G$ is \emph{$b$-positive} if $X_G$ is $b$-positive. It is \emph{nice} if $\lambda\ge\mu$ and $\stc(G;\lambda)>0$ imply $\stc(G;\mu)>0$, and it is \emph{strongly nice} if $\stc(G;\lambda)\le\stc(G;\mu)$ whenever $\lambda\ge\mu$. 

\begin{proposition}[\citeauthor{LLYZ25}]\label{prop:s+:sn+:n+}
Any Schur-positive graph is strongly nice, and any strongly nice graph is nice.
\end{proposition}

The \emph{incomparability graph} $\inc(P)$ of a poset $P$ has vertex set $P$, with two vertices adjacent exactly when they are incomparable. We call $P$ $b$-positive, strongly nice, or nice when $\inc(P)$ has the corresponding property. 

A \emph{lattice} is a poset in which every pair of elements has a meet and a join. A lattice is \emph{distributive} if its meet and join operations satisfy the distributive identities.

For a positive integer $n$, let $\mathbf n$ denote the $n$-element chain $1<2<\dots<n$. For positive integers~$n_1,\,\dots,\,n_r$, let $\mathbf{n_1}\times\cdots\times\mathbf{n_r}$ denote the product of the chains $\mathbf{n_1},\,\dots,\,\mathbf{n_r}$. This product is a graded distributive lattice. In particular, when $n_i=2$ for every $i$, the product $\mathbf 2\times\cdots\times\mathbf 2$ is the Boolean lattice~$B_r$. For $m,n\ge2$, the lattice $\mathbf m\times\mathbf n$ has distinct minimum and maximum elements, which correspond to isolated vertices in the graph $\inc(\mathbf m\times\mathbf n)$. Let $H_m^n$ be the graph obtained from $\inc(\mathbf m\times\mathbf n)$ by removing these two isolated vertices. \Cref{fig:H2} shows the lattice $\mathbf{(m+1)}\times\mathbf 2$ and its reduced incomparability graph~$H_{m+1}^2$, whereas \cref{fig:H3} shows the lattice $\mathbf m\times\mathbf 3$ and its reduced incomparability graph $H_m^3$. By the definition of chromatic symmetric functions,
\begin{equation}\label{inc=H+2}
X_{\inc(\mathbf m\times\mathbf n)}=s_1^2 X_{H_m^n},
\quad\text{where $s_1^2=s_2+s_{11}$.}
\end{equation}
We analyze the product in \cref{inc=H+2} using Pieri's rules; see \citet{Sta99B}.

\begin{proposition}[Pieri's rules]\label{prop:Pieri}
For any partition $\mu$, 
\[
s_\mu s_n=\sum_\lambda s_\lambda
\quad\text{and}\quad
s_\mu s_{1^n}=\sum_\lambda s_\lambda,
\]
where the first sum runs over all partitions $\lambda$ such that $\lambda/\mu$ is a horizontal strip of size $n$, and the second sum runs over all partitions $\lambda$ such that $\lambda/\mu$ is a vertical strip of size $n$.
\end{proposition}

For example, applying Pieri's rule to $s_1s_1$ gives
\[
s_1^2=s_2+s_{11},
\]
corresponding to adding one cell either to the first row or as a new row.

We next record the structure of $H_m^n$. Its vertex and edge sets are
\begin{align*}
V(H_m^n)
&=\{v_{ij}\colon 0\le i\le m-1,\,0\le j\le n-1\}
\setminus \{v_{00},\,v_{(m-1)(n-1)}\}
\quad\text{and}\\
E(H_m^n)
&=\{v_{st}v_{pq}\colon (s-p)(t-q)<0\},
\end{align*}
respectively. It follows that 
\[
\abs{V(H_m^n)}=mn-2
\quad\text{and}\quad
\abs{E(H_m^n)}=\binom{m}{2}\binom{n}{2}.
\]
Its independence number is 
\begin{equation}\label{alphaG}
\alpha\brk1{H_m^n}=m+n-3.
\end{equation}
Indeed, $H_m^n$ has the clique partition
\begin{equation}\label{clique-partition:H3}
\kappa
=\Delta_1
\big/
\Delta_2
\big/
\dotsm
\big/
\Delta_{m+n-3}, 
\quad\text{where $\Delta_k=\{v_{ij}\colon i+j=k\}$,}
\end{equation}
and it contains the stable set 
\[
\brk[c]1{v_{10},\,v_{20},\,\dots,\,v_{(m-1)0},\,v_{(m-1)1},\,\dots,\,v_{(m-1)(n-2)}}
\]
of size $m+n-3$. The clique partition gives $\alpha(H_m^n)\le m+n-3$, while the displayed
stable set gives the reverse inequality. Hence 
\[
\alpha(H_m^n)=m+n-3.
\]
We also need the following description of stable sets in $H_m^n$. 

\begin{lemma}\label{lem:shape:stable}
For any stable set $S$ of the graph $H_m^n$ and for any $0\le q<t\le n-1$, we have
\[
\max\{s\colon v_{sq}\in S\}
\le
\min\{s\colon v_{st}\in S\},
\]
whenever both sides of this inequality are well defined. 
\end{lemma}
\begin{proof}
Since $S$ is stable, any vertices $v_{pq}$ and $v_{st}$ in $S$ are not adjacent in $H_m^n$, i.e., $(p-s)(q-t)\ge 0$. Because $q<t$, it follows that $p\le s$, as required.
\end{proof}

When $v_{ij}$ is drawn at $(i,j)$, \cref{lem:shape:stable} says that every stable set extends monotonically from southwest to northeast. 

For later counts, we refine $\mathcal{SC}(H_m^n;\tau)$ as follows:
\begin{align*}
\mathcal{SC}_{ij}\brk1{H_m^n;\,\tau} 
&=\brk[c]1{(B_1,\, \dots,\, B_{\ell(\tau)})
\in\mathcal{SC}\brk1{H_m^n;\,\tau} 
\colon
\brk1{v_{0(n-1)},\,v_{(m-1)0}}
\in B_i\times B_j}
\quad\text{and}\\
\stc_{ij}\brk1{H_m^n;\,\tau} 
&=\abs{\mathcal{SC}_{ij}\brk1{H_m^n;\,\tau}}.
\end{align*}
We suppress $H_m^n$ and $\tau$ from $\mathcal{SC}_{ij}$ and $\stc_{ij}$ when they are clear from context.

\begin{lemma}\label{lem:split}
For any composition $\tau\vDash mn-2$, we have
$\stc\brk1{H_m^n;\,\tau} 
=2\sum_{1\le i<j\le \ell(\tau)}
\stc_{ij}\brk1{H_m^n;\,\tau}$.
\end{lemma}
\begin{proof}
Since $H_m^n$ has the involutive automorphism $v_{ij}\mapsto v_{(m-1-i)(n-1-j)}$, we obtain $\stc_{ij}=\stc_{ji}$. On the other hand, $\stc_{ii}=0$ since $v_{0(n-1)}v_{(m-1)0}$ is an edge. Summing over the possible ordered block pairs containing $(v_{0(n-1)},v_{(m-1)0})$ gives the formula. 
\end{proof}

We use \cref{lem:split} to count stable compositions of selected types in $H_m^2$ and $H_m^3$. 

For sets $A$ and $B$, write $A-B=A\setminus B$; when $B=\{v\}$, write $A-v$. For a statement $P$, let $\chi(P)$ be $1$ if $P$ is true and $0$ otherwise.

The following roadmap separates the general reduction from the case-specific computations. The case $n\ge4$ requires additional ideas, but the first three steps apply uniformly:

\begin{enumerate}
\item Remove the two isolated vertices of
$\inc(\mathbf m\times\mathbf n)$ and write
\[
X_{\inc(\mathbf m\times\mathbf n)}
=s_1^2X_{H_m^n}.
\]

\item Use Pieri's rules to identify the Schur shapes in
$X_{H_m^n}$ that can contribute to a prescribed Schur coefficient of
$X_{\inc(\mathbf m\times\mathbf n)}$.

\item Use the independence number
\[
\alpha(H_m^n)=m+n-3
\]
to discard every special ribbon tabloid whose content has a part
larger than $m+n-3$.

\item For the surviving contents, use the southwest--northeast
structure in \cref{lem:shape:stable} to enumerate the required stable
compositions. For the target shapes chosen in Sections~3 and~4, this
reduces the coefficient calculations to a short list of explicit
counts.
\end{enumerate}

The first three steps apply to every product of two nontrivial chains. The explicit stable-composition formulas in the fourth step are specific to $n=2,3$; for $n\ge4$, several interacting boundaries arise.

\section{The exact threshold for $\mathbf m\times\mathbf 2$}\label{sec:m2}

To prove \cref{thm:ns:m2}, it is convenient to work with $H_{m+1}^2$. We relabel its vertices by setting 
\[
v_{(k-1)1}=x_{k}
\quad\text{and}\quad
v_{k0}=y_k 
\quad\text{for $1\le k\le m$}, 
\]
and denote the resulting graph by $G_m=G_m(X,Y)$. It has order $2m$. More precisely, let 
\[
X=\{x_1, \dots, x_m\}
\quad\text{and}\quad
Y=\{y_1, \dots, y_m\}.
\]
Then $G_m$ is the bipartite graph with bipartition $(X,Y)$ and edge set 
\[
\bigcup_{i=1}^m \{x_i y_j\colon i\le j\le m\}.
\]
Since $G_{m-1}=H_m^2$, the graph $\inc(\mathbf m\times\mathbf 2)$ is isomorphic to the disjoint union of $G_{m-1}$ and two isolated vertices.
\begin{figure}[h]
\begin{tikzpicture}[scale=0.4]
\draw (0,0)--(-1,1)--(0,2)--(1,1)--cycle;
\draw (0,2)--(1,3)--(2,2)--(1,1)--cycle;
\draw (1,3)--(2,4);
\draw (3,3)--(2,2);
\draw (2,4)--(3,5)--(4,4)--(3,3);
\draw (3,5)--(4,6)--(5,5)--(4,4);
\foreach \x/\y in {0/0, -1/1, 0/2, 1/1, 1/3, 2/2, 3/5, 4/4, 4/6, 5/5} {
  \fill[black] (\x, \y) circle (4pt);}
\foreach \x/\y in {3/4, 2/3, 2.5/3.5} {
  \node[ellipsis] at (\x, \y) {};}

\begin{scope}[xshift=10cm, yshift=3cm, scale=2]
\node[vertex] (a1) at (1,0) {};
\node[vertex] (a2) at (2,0) {};
\node[vertex] (a3) at (3,0) {};
\node[vertex] (am) at (6,0) {};
\node[above] at (a1) {$x_1$};
\node[above] at (a2) {$x_2$};
\node[above] at (a3) {$x_3$};
\node[above] at (am) {$x_m$};
\node[above] at (4.5,0) {$\cdots$};
\node[vertex] (b1) at (1.2, -1) {};
\node[vertex] (b2) at (2.2, -1) {};
\node[vertex] (b3) at (3.2, -1) {};
\node[vertex] (bm) at (6.2, -1) {};
\node[below] at (b1) {$y_1$};
\node[below] at (b2) {$y_2$};
\node[below] at (b3) {$y_3$};
\node[below] at (bm) {$y_m$};
\node[below] at (4.5, -1) {$\cdots$};
\foreach \x in {b1, b2, b3, bm} {
  \draw (a1)--(\x);}
\foreach \x in {b2, b3, bm} {
  \draw (a2)--(\x);}
\foreach \x in {b3, bm} {
  \draw (a3)--(\x);}
\draw (am)--(bm);
\end{scope}
\end{tikzpicture}
\caption{The lattice $\mathbf{(m+1)}\times\mathbf 2$ and its reduced incomparability graph $G_m=H_{m+1}^2$.}
\label{fig:H2}
\end{figure}
The graph $G_m$ is closely related to a half-graph. Erd{\H o}s and Hajnal introduced this term for the bipartite graph $(\{a_1,a_2,\dots\},\{b_1,b_2,\dots\})$ with edges $a_ib_j$ whenever $j>i$. As \citet{Erd84} explained, ``the name half-graph comes from the fact that it can be considered to be half of a complete bipartite graph.''

We begin with a basic counting lemma. 
Whenever a displayed type obtained by subtraction contains a zero part,
that part is omitted. If it contains a negative part, the corresponding
term is interpreted as zero. We also use the convention
$\binom nk=0$ for $k<0$ or $k>n$.

\begin{lemma}\label{lem:half:abc}
For any $1 \le a,b,c\le m$ with $a+b+c=2m$, 
\[
\stc_{12}(G_m;\,abc)
=\begin{dcases*}
\binom{m-1}{c}, & if $a=m$ or $b=m$,\\
\binom{m}{c}, & if $a, b\le m-1$. 
\end{dcases*}
\]
\end{lemma}
\begin{proof}
Let $(A, B, C)\in\mathcal{SC}_{12}$. If $a=m$, then $A=X$, and~$C$ must be a $c$-subset of $Y-y_m$, giving $\binom{m-1}{c}$. By symmetry, the same formula holds when~$b=m$. If $a,b\le m-1$, then $A\subseteq X$ and~$B\subseteq Y$. Each composition is uniquely determined by the $c$-element set of indices appearing in $C$. More explicitly, if~$S=\{r_1<\dots<r_c\}\subseteq[m]$, stability forces 
\[
C\cap Y=\{y_{r_1},\,\dots,\,y_{r_{m-b}}\}
\quad\text{and}\quad
C\cap X=\{x_{r_{m-b+1}},\,\dots,\,x_{r_c}\}.
\]
Conversely, every $c$-subset $S$ gives one such composition. Hence the count is $\binom{m}{c}$. 
\end{proof}

We now derive $\stc(G_m;abc)$ from \cref{lem:half:abc}. We write $k\in[\tau_1,\dots,\tau_\ell]$ to mean that $k$ ranges over the listed entries with multiplicity. For example, $\sum_{k\in[1,2,1]}k=4$. 

\begin{proposition}\label{prop:half:abc}
For any $1 \le a,b,c\le m$ with $a+b+c=2m$, we have 
\[
\stc(G_m;\,abc) 
=\begin{dcases*}
2\brk4{1+\binom{m}{\min(a,b,c)}}, 
& if $\max(a, b, c)=m$, 
\\
2\sum_{k\in[a, b, c]}\binom{m}{k}, 
& if $\max(a, b, c)\le m-1$.
\end{dcases*}
\]
In addition, $\stc(G_m;m^2)=2$.
\end{proposition}
\begin{proof}
By \cref{lem:split},
\[
\stc=2(\stc_{12}+\stc_{13}+\stc_{23}),
\]
and each term is given by \cref{lem:half:abc}. Substitution gives the stated formula. The only two stable compositions of type $m^2$ are $(X,Y)$ and $(Y,X)$.
\end{proof}

We now apply these formulas to specific partitions relevant for the Schur expansion. 

\begin{proposition}\label{prop:s-coeff:half}
For any $m\ge 6$, the graph $G_m$ is not Schur positive. Moreover,
\begin{align*}
[s_{(m-2)^{2}4}]X_{G_m}
&=-m(m-1)(m-5)/3, \\
[s_{(m-1)(m-2)3}]X_{G_m}
&=0, \quad\text{and}\\
[s_{(m-1)^22}]X_{G_m}
&=2m-2.
\end{align*}
For $m\ge7$, one also has
\[
[s_{(m-1)(m-3)4}]X_{G_m}=0.
\]
\end{proposition}
\begin{proof}
Let $\lambda=(m-2)^24\vdash2m$. The six special ribbon tabloids of shape $\lambda$ are shown in \cref{fig:m-2.m-2.4}.
\begin{figure}[h]
\begin{tikzpicture}[scale=0.5, round]
\input H2-sh.tex
\draw[ribbon] (0.5, 0.5)--(1.5, 0.5);
\draw[ribbon] (0.5, 1.5)--(2.5, 1.5)--(2.5, 0.5)--(5.5, 0.5);
\draw[ribbon] (0.5, 2.5)--(3.5, 2.5)--(3.5, 1.5)--(6.5, 1.5)--(6.5, 0.5);
\foreach \x/\y in {5.5/0.5, 1.5/0.5}{\fill[black] (\x, \y) circle (3pt);}
\node at (kappa) {$\kappa_1=m(m-2)2$};

\begin{scope}[xshift=10cm]
\input H2-sh.tex
\draw[ribbon] (0.5, 0.5)--(5.5, 0.5);
\draw[ribbon] (0.5, 1.5)--(2.5, 1.5);
\draw[ribbon] (0.5, 2.5)--(3.5, 2.5)--(3.5, 1.5)--(6.5, 1.5)--(6.5, 0.5);
\foreach \x/\y in {2.5/1.5, 5.5/0.5}{\fill[black] (\x, \y) circle (3pt);}
\node at (kappa) {$\kappa_2=m3(m-3)$};
\end{scope}

\begin{scope}[xshift=20cm]
\input H2-sh.tex
\draw[ribbon] (0.5, 2.5)--(3.5, 2.5)--(3.5, 1.5)--(6.5, 1.5);
\draw[ribbon] (0.5, 1.5)--(2.5, 1.5)--(2.5, 0.5)--(6.5, 0.5);
\draw[ribbon] (0.5, 0.5)--(1.5, 0.5);
\foreach \x/\y in {6.5/1.5, 1.5/0.5}{\fill[black] (\x, \y) circle (3pt);}
\node at (kappa) {$\kappa_3=(m-1)^22$};
\end{scope}

\begin{scope}[yshift=-5cm]
\input H2-sh.tex
\draw[ribbon] (0.5, 2.5)--(3.5, 2.5)--(3.5, 1.5)--(6.5, 1.5);
\draw[ribbon] (0.5, 1.5)--(2.5, 1.5);
\draw[ribbon] (0.5, 0.5)--(6.5, 0.5);
\foreach \x/\y in {2.5/1.5, 6.5/1.5}{\fill[black] (\x, \y) circle (3pt);}
\node at (kappa) {$\kappa_4=(m-1)3(m-2)$};
\end{scope}

\begin{scope}[xshift=10cm, yshift=-5cm]
\input H2-sh.tex
\draw[ribbon] (0.5, 2.5)--(3.5, 2.5);
\draw[ribbon] (0.5, 1.5)--(6.5, 1.5)--(6.5, 0.5);
\draw[ribbon] (0.5, 0.5)--(5.5, 0.5);
\foreach \x/\y in {5.5/0.5, 3.5/2.5}{\fill[black] (\x, \y) circle (3pt);}
\node at (kappa) {$\kappa_5=4(m-1)(m-3)$};
\end{scope}

\begin{scope}[xshift=20cm, yshift=-5cm]
\input H2-sh.tex
\draw[ribbon] (0.5, 0.5)--(6.5, 0.5);
\draw[ribbon] (0.5, 1.5)--(6.5, 1.5);
\draw[ribbon] (0.5, 2.5)--(3.5, 2.5);
\foreach \x/\y in {3.5/2.5, 6.5/1.5}{\fill[black] (\x, \y) circle (3pt);}
\node at (kappa) {$\kappa_6=4(m-2)^2$};
\end{scope}
\end{tikzpicture}
\caption{The set $\mathcal{T}_{(m-2)^24}=\{\kappa_1, \dots, \kappa_6\}$.}
\label{fig:m-2.m-2.4}
\end{figure}
For brevity, we write $N_\kappa=\stc(G_m;\kappa)$. By \cref{thm:scoeff,prop:half:abc}, we obtain
\begin{align*}
[s_\lambda]X_{G_m}
&=-N_{m(m-2)2}
+N_{m(m-3)3}
+N_{(m-1)^22}
-N_{(m-1)(m-2)3}
-N_{(m-1)(m-3)4}
+N_{(m-2)^24}
\\
&=-2\brk3{1+\binom{m}{2}}
+2\brk3{1+\binom{m}{3}}
+2\brk3{2m+\binom{m}{2}}
-2\brk3{m+\binom{m}{2}+\binom{m}{3}}
\\
&\qquad-2\brk3{m+\binom{m}{3}+\binom{m}{4}}
+2\brk3{2\binom{m}{2}+\binom{m}{4}}
\\
&=-\frac{m(m-1)(m-5)}{3},
\end{align*}
which is negative for $m\ge6$. Hence $G_m$ is not Schur positive. The same method gives, for $m\ge6$,
\begin{align*}
[s_{(m-1)(m-2)3}]X_{G_m}
&=N_{m(m-1)1}
-N_{(m-1)^22}
-N_{m(m-3)3}
+N_{(m-1)(m-2)3}
=0, \\
[s_{(m-1)^22}]X_{G_m}
&=N_{m^2}
-N_{m(m-1)1}
-N_{m(m-2)2}
+N_{(m-1)^22}
=2(m-1).
\end{align*}
For $m\ge7$, it also gives
\[
[s_{(m-1)(m-3)4}]X_{G_m}
=N_{m(m-2)2}
-N_{(m-1)(m-2)3}
-N_{m(m-4)4}
+N_{(m-1)(m-3)4}
=0.
\]
This completes the proof.
\end{proof}

We can now prove the first main theorem.

\begin{theorem}\label{thm:ns:m2}
For every positive integer $m$, the lattice $\mathbf m\times\mathbf 2$ is not Schur positive if and only if~$m\ge 8$. Moreover, for every $m\ge 8$,
\[
[s_{(m-2)^24}]X_{\inc(\mathbf m\times\mathbf 2)}
=-m(m-2)(m-7)/3.
\]
\end{theorem}

\begin{proof}
Fix $m\ge8$ and let $\lambda=(m-2)^24$. Let $\mathcal S$ be the set of partitions $\nu$ for which $[s_\lambda]s_\nu s_1^2\ne0$, where $s_1^2=s_2+s_{11}$. By \cref{prop:Pieri}, these partitions are obtained by removing a horizontal or vertical strip of size $2$ from $\lambda$. Thus $\mathcal S=\{\nu_1,\dots,\nu_4\}$, as shown in \cref{fig:H2:nu}. 
\begin{figure}[h]
\begin{tikzpicture}[scale=0.5, round]
\input H2-rm.tex
\draw[gray] (3, 2)--(3, 3);
\draw[border] (0,0)--(6,0)--(6,2)--(4,2)--(4,3)--(0,3)--cycle;
\node[vertex] at (6.5, .5) {};
\node[vertex] at (6.5, 1.5) {};
\node at (3, -.8) {$\nu_1=(m-3)^24$};

\begin{scope}[xshift=8cm]
\input H2-rm.tex
\draw[gray] (6, 0)--(6, 1);
\draw[border] (0,0)--(7,0)--(7,1)--(6,1)--(6,2)--(3,2)--(3,3)--(0,3)--cycle;
\node[vertex] at (3.5, 2.5) {};
\node[vertex] at (6.5, 1.5) {};
\node at (3.5, -.8) {$\nu_2=(m-2)(m-3)3$};
\end{scope}

\begin{scope}[xshift=16cm]
\input H2-rm.tex
\draw[gray] (3, 2)--(3, 3);
\draw[gray] (6, 0)--(6, 1);
\draw[border] (0,0)--(7,0)--(7,1)--(5,1)--(5,2)--(4,2)--(4,3)--(0,3)--cycle;
\node[vertex] at (5.5, 1.5) {};
\node[vertex] at (6.5, 1.5) {};
\node at (3.5, -.8) {$\nu_3=(m-2)(m-4)4$};
\end{scope}

\begin{scope}[xshift=24cm]
\input H2-rm.tex
\draw[gray] (6, 1)--(7, 1);
\draw[gray] (6, 0)--(6, 2);
\draw[border] (0,0)--(7,0)--(7,2)--(2,2)--(2,3)--(0,3)--cycle;
\node[vertex] at (2.5, 2.5) {};
\node[vertex] at (3.5, 2.5) {};
\node at (3.5, -.8) {$\nu_4=(m-2)^22$};
\end{scope}
\end{tikzpicture}
\caption{The contributing partitions $\nu_1,\dots,\nu_4$; shaded circles mark the removed cells.}
\label{fig:H2:nu}
\end{figure}
The partition $\nu_2$ occurs with multiplicity two because the two removed cells form both a horizontal strip and a vertical strip. Let $H=G_{m-1}$. By \cref{prop:s-coeff:half,prop:Pieri}, we then compute
\begin{multline*}
[s_\lambda]X_{\inc(\mathbf m\times\mathbf 2)}
=[s_\lambda]s_1^2X_H
=[s_{\nu_1}]X_H
+2[s_{\nu_2}]X_H
+[s_{\nu_3}]X_H
+[s_{\nu_4}]X_H
\\
=-(m-1)(m-2)(m-6)/3
+2\cdot 0
+0
+(2m-4)
=-m(m-2)(m-7)/3, 
\end{multline*}
which is negative for $m\ge 8$.

The exact computations described in \cref{app:small-cases} show that $\mathbf m\times\mathbf2$ is Schur positive for $m\le7$. This completes the proof.
\end{proof}

\section{The exact threshold for $\mathbf m\times\mathbf 3$}\label{sec:m3}

We now prove \cref{thm:ns:m3}. The graph $\inc(\mathbf m\times\mathbf3)$ is the disjoint union of $H_m^3$ and two isolated vertices; see \cref{fig:H3}.

\begin{figure}[h]
\begin{tikzpicture}[scale=0.45]
\draw (0, 0)--(-2, 2)--(3, 7)--(5, 5)--cycle;
\draw (-1, 3)--(1, 1);
\draw (0, 4)--(2, 2);
\draw (2, 6)--(4, 4);
\draw (-1, 1)--(4, 6);
\foreach \x/\y in {0/0, -1/1, -2/2, 1/1, 0/2, -1/3, 2/2, 1/3, 0/4, 4/4, 3/5, 2/6, 5/5, 4/6, 3/7} {
  \fill[black] (\x, \y) circle (4pt);}
\foreach \x/\y in {2/3, 2.5/3.5, 3/4, 1/4, 1.5/4.5, 2/5} {
  \node[ellipsis] at (\x, \y) {};}
\begin{scope}[yshift=5cm, xshift=9cm, scale=2]
\node[vertex] (x0) at (1,0) {};
\node[vertex] (x1) at (2,0) {};
\node[vertex] (x2) at (3,0) {};
\node[vertex] (xm) at (6,0) {};
\node[above] at (x0) {$x_0$};
\node[above] at (x1) {$x_1$};
\node[above] at (x2) {$x_2$};
\node[above] at (xm) {$x_{m-2}$};
\node[above] at (4.5, 0) {$\cdots$};
\node[vertex] (y0) at (0.5, -1) {};
\node[vertex] (y1) at (1.5,- 1) {};
\node[vertex] (y2) at (2.5, -1) {};
\node[vertex] (ym) at (5.5, -1) {};
\node[vertex] (yn) at (6.5, -1) {};
\node[left] at (y0) {$y_0$};
\node[right] at (yn) {$y_{m-1}$};
\node[vertex] (z1) at (1.25, -2) {};
\node[vertex] (z2) at (2.25, -2) {};
\node[vertex] (z3) at (3.25, -2) {};
\node[vertex] (zm) at (6.25, -2) {};
\node[below] at (z1) {$z_1$};
\node[below] at (z2) {$z_2$};
\node[below] at (z3) {$z_3$};
\node[below] at (zm) {$z_{m-1}$};
\node[below] at (4.5, -2) {$\cdots$};
\foreach \x in {y1, y2, ym, yn, z1, z2, z3, zm} {
  \draw (x0)--(\x);}
\foreach \x in {x1, x2, xm, y0, y1, y2, ym} {
  \draw (zm)--(\x);}
\foreach \x in {y2, ym, yn, z2, z3} {
  \draw (x1)--(\x);}
\foreach \x in {ym, yn, z3} {
  \draw (x2)--(\x);}
\foreach \x in {z1, z2, z3} {
  \draw (y0)--(\x);}
\draw (xm)--(yn);
\draw (y1)--(z2);
\draw (y1)--(z3);
\draw (y2)--(z3);
\end{scope}
\end{tikzpicture}
\caption{The lattice $\mathbf m\times\mathbf 3$ and the graph $H_m^3$.}
\label{fig:H3}
\end{figure}
For convenience, we relabel the vertices by setting $v_{k2}=x_k$ for~$0\le k\le m-2$, $v_{k1}=y_k$ for~$0\le k\le m-1$, and $v_{k0}=z_k$ for~$1\le k\le m-1$, and denote the relabeled graph by $H_m$. The graph has order $3m-2$. More precisely, let 
\[
X=\{x_0,\dots,x_{m-2}\},
\quad
Y=\{y_0,\dots,y_{m-1}\}
\quad\text{and}\quad
Z=\{z_1,\dots,z_{m-1}\}.
\]  
Then $H_m$ has vertex set $X\sqcup Y\sqcup Z$ and edge set 
\[
\{x_i y_j,\,y_i z_j,\,x_i z_j\colon 0\le i<j\le m-1\}.
\]  
For any pair $(W,w)\in\{(X,x),(Y,y),(Z,z)\}$ and integers
$i,j$, write
\[
W_{ij}:=\{w_r:w_r\in W,\ i\le r\le j\}.
\]
In particular, $W_{ij}=\emptyset$ whenever $i>j$, $X=X_{0(m-2)}$, $Y=Y_{0(m-1)}$ and $Z=Z_{1(m-1)}$.	

We first count stable compositions of $H_m$ whose type has length $3$. 

\begin{proposition}\label{prop:len.tau=3}
For $m\ge 3$, we have $\stc\brk1{H_m;\,m^2(m-2)}=6$ and $\stc\brk1{H_m;\,m(m-1)^2}=12$.
\end{proposition}
\begin{proof}
Let $\tau=m^2(m-2)$ and $(A, B, C)\in\mathcal{SC}_{12}(H_m;\tau)$. By \cref{lem:shape:stable}, we have
\begin{align}\label{A=yX:m}
A&=\{y_0\}\cup X
\quad\text{and}\\
\label{B=Zy:m-1}
B&=Z\cup\{y_{m-1}\}.
\end{align}
Thus $\stc_{12}(\tau)=1$, and the same argument gives $\stc_{13}(\tau)=\stc_{23}(\tau)=1$. Hence \cref{lem:split} yields
\[
\stc(\tau)
=2\brk1{\stc_{12}(\tau)+\stc_{13}(\tau)+\stc_{23}(\tau)}=6.
\] 
For $\eta=m(m-1)^2$, a similar enumeration gives
\[
\stc(\eta)=2\brk1{\stc_{12}(\eta)+\stc_{13}(\eta)+\stc_{23}(\eta)}=2\cdot(2+2+2)=12.
\] 
This completes the proof.
\end{proof}

We next consider length-$4$ types with exactly two parts equal to $m$.

\begin{proposition}\label{prop:m.m.c.d}
For any $1\le c,d\le m-3$ such that $c+d=m-2$, 
\[
\stc\brk1{H_m;\,m^2cd}
=6\binom{m-2}{c}+12.
\]
\end{proposition}
\begin{proof}
Let $\tau=m^2cd$ and $\pi=(A,B,C,D)\in\mathcal{SC}$.

If $\pi\in \mathcal{SC}_{12}$, then \cref{lem:shape:stable} implies \cref{A=yX:m,B=Zy:m-1}, and $C\sqcup D=Y_{1(m-2)}$. Thus $\stc_{12}=\binom{m-2}{c}$ for the choices of $C$. If $\pi\in \mathcal{SC}_{13}$, then \cref{lem:shape:stable} implies \cref{A=yX:m}, and $z_1,y_{m-1}\in B$. Then the restriction of $\pi$ to the subgraph
\begin{equation}\label{G-X-z-y:m}
H_m-A-\{z_1,\,y_{m-1}\}\cong 
G_{m-2}\brk1{Y_{1(m-2)},\,Z_{2(m-1)}}
\end{equation}
has type $(m-2)cd$, with $(y_1, z_{m-1})\in(B\cup D)\times C$. By \cref{lem:half:abc}, 
\begin{equation}\label{pf:s13:mmcd}
\stc_{13}
=\stc_{12}\brk1{G_{m-2};\,(m-2)cd} 
+\stc_{12}\brk1{G_{m-2};\,dc(m-2)}
=\binom{m-3}{d}+1.
\end{equation}
If $\pi\in\mathcal{SC}_{34}$, then by \cref{lem:shape:stable},
\[
\{A,B\}
=\brk[c]1{
Y_{0c}\cup X_{c(m-2)},\ Z_{1(c+1)}
\cup Y_{(c+1)(m-1)}
}
\quad\text{and}\quad
(C,D)=\brk1{X_{0(c-1)},\,Z_{(c+2)(m-1)}}.
\]
Thus $\stc_{34}=2$, corresponding to the two possible ordered pairs $(A,B)$.

By exchanging $c$ and $d$ in \cref{pf:s13:mmcd}, we obtain a formula for $\stc_{14}$. By symmetry, we have $\stc_{23}=\stc_{13}$ and $\stc_{24}=\stc_{14}$. Hence by \cref{lem:split}, 
\[
\stc
=2\sum_{1\le i<j\le 4}\stc_{ij}
=2\brk4{
\binom{m-2}{c}
+
2\brk3{\binom{m-3}{d}+1
+\binom{m-3}{c}+1}
+2
},
\]
which simplifies to the stated formula.
\end{proof}

We next consider length-$4$ types with exactly one part equal to $m$. In the formulas below, a zero part is omitted, a term containing a negative part is interpreted as zero, and 
\[
\binom nk=0\quad\text{for $k<0$ or $k>n$}.
\]

\begin{lemma}\label{lem:s1j:mbcd}
Let $m\ge 5$ and $1\le b,c,d\le m-1$ such that $b+c+d=2m-2$. Let $s=\max(b,c,d)$ and $t=\min(b,c,d)$. Then the sum $\sum_{j=2}^4 \stc_{1j}(H_m;\,mbcd)$ equals
\[
\begin{dcases*}
6\binom{m-1}{t}
+4
+2m\chi(t\ge 2),
& if $s=m-1$, 
\\
4\brk4{\binom{m}{t}+\binom{m-1}{t-1}+\binom{m-2}{t}}
+3m
-2
+m\chi(t=2)
+m^2\chi(t\ge 3),
& if $s=m-2$,
\\
2\sum_{k\in[b,c,d]}\brk4{\binom{m}{k}+2\binom{m-1}{k}+\binom{m-2}{k}},
& otherwise.
\end{dcases*}
\]
As above, $\chi(P)$ denotes the indicator of the statement $P$. 
\end{lemma}
\begin{proof}
Let $\pi=(A,B,C,D)\in\mathcal{SC}_{12}\cup\mathcal{SC}_{13}\cup\mathcal{SC}_{14}$. By \cref{lem:shape:stable}, we have \cref{A=yX:m,G-X-z-y:m}. Distributing the isolated vertices $z_1$ and $y_{m-1}$ among the remaining blocks gives
\begin{multline*}
\sum_{j=2}^4 \stc_{1j}
=\sum_{
\{\epsilon_b,\,\epsilon_c,\,\epsilon_d\}=\{0,0,2\}\text{ as multisets}}
\stc\brk1{G_{m-2};\,(b-\epsilon_b)(c-\epsilon_c)(d-\epsilon_d)}
\\
+\sum_{
\{\epsilon_b,\,\epsilon_c,\,\epsilon_d\}=\{0,1,1\}\text{ as multisets}}
2\cdotp
\stc\brk1{G_{m-2};\,(b-\epsilon_b)(c-\epsilon_c)(d-\epsilon_d)},
\end{multline*}
where the factor $2$ comes from the distinctness of $z_1$ and $y_{m-1}$. By \cref{prop:half:abc}, substituting the preceding expressions into the formula above and simplifying gives the stated formula.
\end{proof}

We now specialize to types containing both $m$ and $m-1$. 

\begin{proposition}\label{prop:m.m-1.c.d}
Let $m\ge 4$, $1\le c,d\le m-2$, and $c+d=m-1$. Then
\[
\stc\brk1{H_m;\,m(m-1)cd}
=18\binom{m-1}{c}+34+(8m+2)\chi(\min(c, d)\ge 2).
\]
\end{proposition}
\begin{proof}
The case $m=4$ is verified directly. Hence assume $m\ge 5$, so that \cref{lem:s1j:mbcd} applies.

Write $a=m$, $b=m-1$, and $\tau=abcd$. By \cref{lem:split,lem:s1j:mbcd}, it remains to compute $\stc_{23}$, $\stc_{24}$, and $\stc_{34}$. Let $\pi=(A,B,C,D)\in\mathcal{SC}$. 

For any $\pi\in\mathcal{SC}_{23}$ with $B\ne X$, \cref{lem:shape:stable} and $\abs{B}=m-1$ imply
\[
B=\{y_0\}\cup(X-x_i)
\quad\text{for a unique $i\in[m-2]$}.
\]
We decompose $\mathcal{SC}_{23}=\mathcal N_1\sqcup\mathcal N_2\sqcup\mathcal N_3$, where
\begin{align*}
\mathcal N_1
&=\{\pi=(A,B,C,D)\in\mathcal{SC}_{23}\colon B=X\}, \\
\mathcal N_2
&=\{\pi=(A,B,C,D)\in\mathcal{SC}_{23}\colon B\ne X \text{ and }
x_i\in A\}, \quad\text{and}\\
\mathcal N_3
&=\{\pi=(A,B,C,D)\in\mathcal{SC}_{23}\colon B\ne X
\text{ and }x_i\in D\}.
\end{align*}
Let $N_i=\abs{\mathcal N_i}$ for $i=1,2,3$; we compute these three quantities in turn. 
\begin{enumerate}
\item
Let $\pi\in\mathcal N_1$. Then the restriction of $\pi$ to the graph
\begin{equation}\label{G-X}
H_m-X\cong
G_{m-1}(Y-y_{m-1},\,Z)\cup\{y_{m-1}\}
\end{equation}
has type $acd$, with $(y_0, z_{m-1})\in (A\cup D)\times C$. By \cref{lem:half:abc},
\[
N_1
=\stc_{12}\brk1{G_{m-1};\,(a-1)cd}
+\stc_{12}\brk1{G_{m-1};\,dc(a-1)}
=\binom{m-2}{d}+1.
\]
\item
Let $\pi\in\mathcal N_2$. Then 
\begin{equation}\label{B=yX-xi:m-1}
B=\{y_0\}\cup(X-x_i),
\quad i\in[m-2].
\end{equation}
By \cref{lem:shape:stable}, necessarily $i=m-2$, and
\[
A=Z_{1k}\cup Y_{k(m-2)}\cup\{x_{m-2}\}
\quad\text{for some $k\in[m-2]$.}
\]
If $k=1$, then $D$ is a $d$-subset of $Z_{2(m-2)}\cup\{y_{m-1}\}$. If $k\ge 2$, then
\[
(k, C, D)\in\brk[c]1{
\brk1{d, \,
Z_{(d+1)(m-1)},\,
Y_{1(d-1)}\cup\{y_{m-1}\}
}, \
\brk1{d+1, \,
Z_{(d+2)(m-1)}\cup\{y_{m-1}\},\,
Y_{1d}
}}.
\]
Since $2\le k\le m-2$, one of the two possibilities requires $d\ge 2$ and the other requires $c\ge 2$; thus
\[
N_2=\binom{m-2}{d}+1+\chi(t\ge 2), 
\quad\text{where $t=\min(c, d)$.}
\]
\item
Let $\pi\in\mathcal N_3$. Then \cref{B=yX-xi:m-1} holds. By \cref{lem:shape:stable}, we find
\[
A=Z_{1k}\cup Y_{k(m-1)}
\quad\text{for some $1\le k\le \min(m-2,\,i+1)$}.
\]
It follows that 
\[
(C,D)=\brk1{C'\cup Z_{(i+1)(m-1)},\
D'\cup Y_{1(k-1)}\cup\{x_i\}},
\quad\text{where $C'\sqcup D'=Z_{(k+1)i}$.}
\]
If $k=1$, then $D'$ is a $(d-1)$-subset of $Z_{2i}$. If $k\ge 2$, then $D'=\emptyset$ since~$D$ is stable; thus 
\[
(k, C,D)=\brk1{d, \, 
Z_{(d+1)(m-1)},\, 
Y_{1(d-1)}\cup\{x_i\}
}.
\]
Since $d-1\le i\le m-2$, we obtain
\[
N_3
=\sum_{i=1}^{m-2}\binom{i-1}{d-1}
+(m-d)\chi(d\ge 2)
=\binom{m-2}{d}+(c+1)\chi(d\ge 2).
\]
\end{enumerate}
Combining the three cases gives
\begin{equation}\label{pf:s23:m.m-1.c.d}
\stc_{23}
=N_1+N_2+N_3
=3\binom{m-2}{d}+2+\chi(t\ge 2)+(c+1)\chi(d\ge 2).
\end{equation}

We next compute $\stc_{34}$. Let $\pi\in \mathcal{SC}_{34}$. By \cref{lem:shape:stable}, 
\begin{equation}\label{A:mbcd}
A=Z_{1k}\cup Y_{kj}\cup X_{j(m-2)},
\end{equation}
for some $0\le k\le m-2$ and $1\le j\le m-1$ such that $k\le j$. We divide the count according to $k$ and~$j$, writing $N_1',\dots,N_4'$ for the four contributions.
\begin{enumerate}[
label=Case \arabic*), 
leftmargin=*,
labelwidth=\casewidth,
itemsep=5pt,
]
\item
If $k=0$ and $j=m-1$, then $A=Y$ and
\begin{equation}\label{G-Y}
H_m-Y\cong G_{m-1}(X,Z).
\end{equation}
Thus, by \cref{lem:half:abc}, $N_1'=\stc_{12}\brk1{G_{m-1};\,cd(m-1)}=1$.
\item
If $k=0$ and $j\le m-2$, then $A=Y_{0j}\cup X_{j(m-2)}$. By \cref{lem:shape:stable} and the condition $\abs{B}=m-1$, we infer that $X\cap B=\emptyset$. It follows that $j=c$ and $C=X_{0(c-1)}$. If $d=1$, then 
\[
\pi=\brk1{
Y_{0(m-2)}\cup\{x_{m-2}\},\
Z_{1(m-2)}\cup\{y_{m-1}\},\
X_{0(m-3)},\
\{z_{m-1}\}}
\quad\text{and}\quad
N_2'=1.
\]
If $d\ge 2$, then $c\le m-3$ and $Y_{(c+1)(m-2)}\subseteq B$. It follows that $Z_{(c+2)(m-1)}\subseteq D$ and the remaining vertex of $D$ can be chosen from $Z_{1(c+1)}\cup\{y_{m-1}\}$. Therefore, 
\begin{equation}\label{pf:s34:k=0.j<m-1:m.m-1.c.d}
N_2'=1+(c+1)\chi(d\ge 2).
\end{equation}
\item
If $k\ge1$ and $j=m-1$, exchanging $c$ and $d$ in \cref{pf:s34:k=0.j<m-1:m.m-1.c.d} gives
\[
N_3'=1+(d+1)\chi(c\ge 2).
\]
\item
Suppose that $1\le k\le j\le m-2$. By \cref{lem:shape:stable,A:mbcd}, the block $B$ is contained in at most two of the sets $X$, $Y$ and~$Z$. If $B$ is contained in a single one of the layers $X$, $Y$, and~$Z$, then that layer must be $Y$, and
\[
(B,C,D)=\brk1{Y-y_j,\,X_{0(j-1)},\,Z_{(k+1)(m-1)}}
\quad\text{and}\quad
j=k=c.
\] 
Otherwise, $B$ is contained in the union of two of the sets $X$, $Y$ and $Z$. Since $z_1,x_{m-2}\in A$, we find $B\not\subseteq X\cup Z$. If $B\subset X\cup Y$, then $j=m-2$, $k=c+1$, and $d\ge 2$. In this case, 
\[
(B,C,D)=\brk1{
Y_{0(k-1)}\cup X_{(k-1)(m-3)},\
X_{0(k-2)},\
Z_{(k+1)(m-1)}\cup\{y_{m-1}\}}.
\]
By symmetry, there exists another possibility for $\pi$ with $B\subset Y\cup Z$ and $c\ge 2$. Therefore, 
\[
N_4'=1+\chi(c\ge 2)
+\chi(d\ge 2)
=2+\chi(t\ge 2).
\]
\end{enumerate}
Therefore, $\stc_{34}=N_1'+\dots+N_4'=7+m\chi(t\ge2)$.
Exchanging $c$ and $d$ in \cref{pf:s23:m.m-1.c.d} gives $\stc_{24}$. By \cref{lem:split,lem:s1j:mbcd}, 
\begin{multline*}
\stc
=2(\stc_{12}+\stc_{13}+\stc_{14})
+2\stc_{23}+2\stc_{24}+2\stc_{34}
\\
=2\brk4{6\binom{m-1}{t}
+4
+2m\chi(t\ge 2)}
+2\brk4{
3\binom{m-2}{d}+2
+\chi(t\ge 2)+(c+1)\chi(d\ge 2)
}
\\
+2\brk4{3\binom{m-2}{c}+2
+\chi(t\ge 2)+(d+1)\chi(c\ge 2)}
+2\brk1{7+m\chi(t\ge 2)},
\end{multline*}
which simplifies to the desired formula.
\end{proof}

We next derive the corresponding formulas when $b\le m-2$.

\begin{lemma}\label{lem:s23:mbcd}
Let $m\ge 5$ and $2\le b,c,d\le m-2$ such that $b+c+d=2m-2$. Let $N=\stc_{23}(H_m;\,mbcd)$. Then we have the following.
\begin{enumerate}
\item
If $b=m-2$, then
\[
N=\begin{dcases*}
3m+7, & if $c=2$, 
\\
3\binom{m-1}{d}
+\binom{m-3}{d}
+\binom{c+2}{2}
+c+3, & if $3\le c\le m-3$,
\\
2m^2-6m+7, & if $c=m-2$.
\end{dcases*}
\]
\item
If $b,c\le m-3$, then
\[
N=4+
\binom{m-1}{d}
+2\binom{m-2}{d}
+\binom{m-3}{d}
+\sum_{k\in[b,c]}
\brk4{\binom{k+2}{m-d}
+\binom{k}{m-d-1}}
-\chi(d=m-2).
\]
\end{enumerate}
\end{lemma}
\begin{proof}
Let $\pi=(A,B,C,D)\in\mathcal{SC}_{23}$. By \cref{lem:shape:stable}, we have \cref{A:mbcd}. Accordingly, we can decompose $X$, $Y$ and $Z$ as 
\begin{align*}
X&=X_{0(k-2)}\sqcup X_{(k-1)(j-1)}\sqcup X_{j(m-2)}, \\
Y&=Y_{1(k-1)}\sqcup Y_{kj}\sqcup Y_{(j+1)(m-2)}, \quad\text{and}\\
Z&=Z_{1k}\sqcup Z_{(k+1)(j+1)}\sqcup Z_{(j+2)(m-1)}, 
\end{align*}
where $Z_{(k+1)(j+1)},\,Y_{kj},\,X_{(k-1)(j-1)}\ne \emptyset$. We first establish the containment relations between these nine sets and the blocks of $\pi$, as illustrated in \cref{fig:decomp:XYZ}.

\begin{figure}[h]
\begin{tikzpicture}
\coordinate (z1) at (0, 0);
\node [above] at (z1.north) {$z_1$};
\coordinate (zk) at (3, 0);
\node [above] at (zk.north) {$z_k$};
\coordinate (zk1) at (4, 0);
\node [above] at (zk1.north) {$z_{k+1}$};
\coordinate (zj1) at (7, 0);
\node [above] at (zj1.north) {$z_{j+1}$};
\coordinate (zj2) at (8, 0);
\node [above] at (zj2.north) {$z_{j+2}$};
\coordinate (zm1) at (11, 0);
\node [above] at (zm1.north) {$z_{m-1}$};

\node[circle, fill=black, inner sep=1pt] (y0) at (-1.5, 1) {};
\node [above] at (y0.north) {$y_0$};
\coordinate (y1) at (-0.5, 1);
\node [above] at (y1.north) {$y_1$};
\coordinate (yk1) at (2.5, 1);
\node [above] at (yk1.north) {$y_{k-1}$};
\coordinate (yk) at (3.5, 1);
\node [above] at (yk.north) {$y_k$};
\coordinate (yj) at (6.5, 1);
\node [above] at (yj.north) {$y_j$};
\coordinate (yj1) at (7.5, 1);
\node [above] at (yj1.north) {$y_{j+1}$};
\coordinate (ym2) at (10.5, 1);
\node [above] at (ym2.north) {$y_{m-2}$};
\node[circle, fill=black, inner sep=1pt] (ym1) at (11.5, 1) {};
\node [above] at (ym1.north) {$y_{m-1}$};

\coordinate (x0) at (-1, 2);
\node [above] at (x0.north) {$x_0$};
\coordinate (xk2) at (2, 2);
\node [above] at (xk2.north) {$x_{k-2}$};
\coordinate (xk1) at (3, 2);
\node [above] at (xk1.north) {$x_{k-1}$};
\coordinate (xj1) at (6, 2);
\node [above] at (xj1.north) {$x_{j-1}$};
\coordinate (xj) at (7, 2);
\node [above] at (xj.north) {$x_j$};
\coordinate (xm2) at (10, 2);
\node [above] at (xm2.north) {$x_{m-2}$};

\draw (z1)--(zk);
\draw (zk1)--(zj1);
\draw (zj2)--(zm1);

\draw (y1)--(yk1);
\draw (yk)--(yj);
\draw (yj1)--(ym2);

\draw (x0)--(xk2);
\draw (xk1)--(xj1);
\draw (xj)--(xm2);
\draw[ultra thick] (z1) to node [above] {$A$} (zk);
\draw[ultra thick] (yk) to node [above] {$A$} (yj);
\draw[ultra thick] (xj) to node [above] {$A$} (xm2);
\draw[ultra thick, dashed] (x0) to node [above] {$B$} (xk2);
\draw[ultra thick, densely dotted, line cap=round] (zj2) to node [above] {$C$} (zm1);
\draw[ultra thick, dash pattern=on 4pt off 2pt on 1pt off 2pt] (y1) to node [above] {$D$} (yk1);
\draw[ultra thick, dash pattern=on 4pt off 2pt on 1pt off 2pt] (yj1) to node [above] {$D$} (ym2);
\end{tikzpicture}
\caption{The decompositions of $X$, $Y$ and $Z$. Thick solid, dashed, dotted, and dash-dotted segments indicate the parts lying in $A$, $B$, $C$, and $D$, respectively.}
\label{fig:decomp:XYZ}
\end{figure}

In fact, \cref{A:mbcd} allows us to locate the subsets $A\cap X$, $A\cap Y$ and $A\cap Z$. By \cref{lem:shape:stable} and the condition $(x_0,z_{m-1})\in B\times C$, we deduce that
\[
B\subseteq\{y_0\}\cup X_{0(j-1)}, \quad
C\subseteq Z_{(k+1)(m-1)}\cup\{y_{m-1}\}, 
\quad\text{and}\quad
Y_{1(k-1)}\cup Y_{(j+1)(m-2)}\subseteq D.
\]
We claim that
\[
X_{0(k-2)}\subseteq B
\quad\text{and}\quad
Z_{(j+2)(m-1)}\subseteq C. 
\]
In fact, $X_{0(k-2)}\cap A=\emptyset$ since $z_k\in A$; and $X_{0(k-2)}\cap C=\emptyset$ since $z_{m-1}\in C$. If $x_i\in D$ for some~$1\le i\le k-2$, then $k\ge 3$; this would imply that $y_{k-1}\in D$, and $D$ would not be stable. This contradiction proves the first statement in the claim. By symmetry, we obtain the other statement. 

These are the containment relations shown in \cref{fig:decomp:XYZ}. We now split into four cases according to whether $X\cap D$ and $Z\cap D$ are empty. Let $N_1,\dots,N_4$ denote the corresponding numbers of compositions. 
\begin{enumerate}[
label=Case \arabic*), 
leftmargin=*,
labelwidth=\casewidth,
itemsep=5pt,
]
\item\label[itm]{itm:DX<>0:DZ<>0}
Suppose that $X\cap D\ne\emptyset$ and $Z\cap D\ne\emptyset$. The condition $X\cap D\ne\emptyset$ implies the following:
\begin{itemize}[topsep=3pt]
\item
The set~$Y_{(j+1)(m-2)}$ is empty, i.e., $j\in\{m-2,\,m-1\}$; 
\item
The set~$Z_{(j+2)(m-1)}$ is empty since $j\ge m-2$; and
\item
The vertex~$y_{m-1}$ does not belong to $D$ since it is adjacent to every vertex of~$X$.
\end{itemize}
By symmetry, the condition $Z\cap D\ne \emptyset$ implies
\[
k\in\{0,1\}, \quad
X_{0(k-2)}=Y_{1(k-1)}=\emptyset, 
\quad\text{and}\quad
y_0\not\in D.
\]
Since two parts in the decomposition of $Y$ are empty, we find $Y_{1(m-2)}\subset A$. We split according to $A$. Recall from the clique partition in \cref{clique-partition:H3} that 
\begin{equation}\label{Delta1:Deltam:H3}
\begin{aligned}
\Delta_1&=\{y_0,z_1\},
&\Delta_k&=\{x_{k-2},y_{k-1},z_k\}\quad(2\le k\le m-1),\\
\Delta_m&=\{x_{m-2},y_{m-1}\}.
\end{aligned}
\end{equation}
Thus $A-Y_{1(m-2)}$ contains one vertex from each of $\Delta_1$ and $\Delta_m$. Let $N_{11},\dots,N_{14}$ count the four corresponding cases. Then 
\[
N_1=N_{11}+\dots+N_{14}.
\] 
\begin{enumerate}[
label=Case 1-\arabic*), 
topsep=5pt, 
leftmargin=\casewidth,
itemsep=5pt
]
\item\label[itm]{itm:DX<>0:DZ<>0:yy}
If $A=Y$, then \cref{G-Y} holds and the restriction of $\pi$ to 
$H_m-Y$ has type $bcd$ with $(x_0,\,z_{m-1})\in B\times C$. By \cref{lem:half:abc}, $N_{11}=\binom{m-1}{d}$.
\item\label[itm]{itm:DX<>0:DZ<>0:yx}
If $A=Y_{0(m-2)}\cup\{x_{m-2}\}$, then $y_{m-1},\,z_{m-2}\in C$ since $X\cap D\ne \emptyset$. On the other hand, since~$X$ contains a vertex in $A$ and a vertex in $D$, we find $N_{12}=0$ if $b\ge m-2$. When~$b\le m-3$, the restriction of $\pi$ to the graph
\[
H_m-Y-\{x_{m-2},\, z_{m-1}\}
\cong
G_{m-2}\brk1{X_{0(m-3)},\ Z_{1(m-2)}}
\]
has type $b(c-2)d$, with $(x_0,\,z_{m-2})\in B\times C$. By \cref{lem:half:abc},
\begin{equation}\label{pf:N12}
N_{12}=\binom{m-2}{d}\chi(b\le m-3).
\end{equation}

\item\label[itm]{itm:DX<>0:DZ<>0:yz}
When $A=\{z_1\}\cup Y_{1(m-1)}$, 
exchanging $b$ and $c$ in \cref{pf:N12} yields 
\[
N_{13}=\binom{m-2}{d}\chi(c\le m-3).
\]
\item\label[itm]{itm:DX<>0:DZ<>0:xz}
If $A=\{z_1\}\cup Y_{1(m-2)}\cup\{x_{m-2}\}$, then, in addition to $x_0\in B$ and $z_{m-1}\in C$, the argument in \cref{itm:DX<>0:DZ<>0:yx} shows that $y_0,\,x_1\in B$ and $y_{m-1},\,z_{m-2}\in C$. Now, the restriction of $\pi$ to 
\[
H_m-Y-\{x_0,\,z_1, \,z_{m-1}, \,x_{m-2}\}
\cong
G_{m-3}\brk1{X_{1(m-3)},\,Z_{2(m-2)}}
\]
has type $(b-2)(c-2)d$ with $(x_1,z_{m-2})\in B\times C$. By \cref{lem:half:abc}, $N_{14}=\binom{m-3}{d}$.
\end{enumerate}

\item\label[itm]{itm:DX<>0:DZ=0}
If $X\cap D\ne \emptyset$ and $Z\cap D=\emptyset$, then 
\[
j\in\{m-2,\,m-1\}, \quad
Y_{(j+1)(m-2)}
=Z_{(j+2)(m-1)}
=\emptyset
\quad\text{and}\quad
y_{m-1}\not\in D,
\] 
by the arguments at the beginning of \cref{itm:DX<>0:DZ<>0}. Let $N_{21}$ and $N_{22}$ denote the numbers of such compositions for $j=m-1$ and $j=m-2$, respectively. Then $N_2=N_{21}+N_{22}$.
\begin{enumerate}[
label=Case 2-\arabic*), 
topsep=5pt, 
leftmargin=\casewidth,
itemsep=5pt
]
\item\label[itm]{itm:DX<>0:DZ=0:j=m-1}
If $j=m-1$, i.e., if $y_{m-1}\in A$, then $C=Z_{(k+1)(m-1)}$ and
\begin{equation}\label{pf:c=m-1-k}
c=m-1-k.
\end{equation}
We separate the possible values of $c$.
\begin{enumerate}[
label=Case 2-1-\arabic*), 
topsep=5pt, 
leftmargin=\casewidth,
itemsep=5pt
]
\item\label[itm]{itm:DX<>0:DZ=0:j=m-1:k<=1}
If $c=m-2$, then $b+d=m$, $k=1$ by \cref{pf:c=m-1-k}. By definition, we have
\[
Y_{1(m-1)}\subset A
\quad\text{and}\quad
X\subseteq B\cup D\subseteq X\cup\{y_0\}.
\]
Since $D$ is a $d$-subset of the stable set $(B\cup D)-x_0$, we find $N_{21}=\binom{m-1}{d}$.

\item\label[itm]{itm:DX<>0:DZ=0:j=m-1:k>=2}
If $c\le m-3$, then $k\ge 2$ by \cref{pf:c=m-1-k}. By definition, 
\[
Y_{1(k-1)}\ne \emptyset
\quad\text{and}\quad
B\cup (D-Y_{1(k-1)})
=\{y_0\}\cup X.
\]
Then $D-Y_{1(k-1)}$ is a $(d-k+1)$-subset of $\{y_0\}\cup X_{(k-1)(m-2)}$, and $N_{21}=\binom{c+2}{m-b}$.
\end{enumerate}
In summary, 
\[
N_{21}=\binom{c+1+\chi(c\le m-3)}{m-b}.
\]
\item\label[itm]{itm:DX<>0:DZ=0:j=m-2}
If $j=m-2$, then $x_{m-2}\in A$, $C=Z_{(k+1)(m-1)}\cup\{y_{m-1}\}$, and $c=m-k$. Since $c\le m-2$, i.e., $k\ge 2$, the set $D-Y_{1(k-1)}$ consists of a single vertex in $X_{(k-1)(m-3)}$ if $b=m-2$, and is a $(d-k+1)$-subset of $\{y_0\}\cup X_{(k-1)(m-3)}$ if $b\le m-3$; thus 
\[
N_{22}=\binom{c}{m-1-d}-\chi(b=m-2).
\]
\end{enumerate}
In summary, 
\begin{equation}\label{N2:mbcd}
N_2=N_{21}+N_{22}
=\binom{c+2-\chi(c=m-2)}{m-d-\chi(c=m-2)}
+\binom{c}{m-d-1}
-\chi(b=m-2).
\end{equation}

\item
If $X\cap D=\emptyset$ and $Z\cap D\ne \emptyset$, then exchanging~$b$ and $c$ in \cref{N2:mbcd} yields $N_3$. 

\item
Suppose that $X\cap D=Z\cap D=\emptyset$. If $y_0\in A$, then $C=Z$ and $c=m-1$, a contradiction. Thus $y_0\in B\cup D$. By symmetry, we have $y_{m-1}\in C\cup D$. We split according to the ordered block pair containing the vertex pair $p=(y_0,y_{m-1})$. Let $N_{41},\dots,N_{44}$ be the corresponding numbers of compositions.
\begin{enumerate}[
label=Case 4-\arabic*), 
topsep=5pt, 
leftmargin=\casewidth,
itemsep=5pt
]
\item
If $p\in B\times C$, then 
\[
(B,C)=\brk1{
\{y_0\}\cup X_{0(j-1)},\
Z_{(k+1)(m-1)}\cup\{y_{m-1}\}
}
\quad\text{and}\quad
(j,k)=(b-1,\,m-c).
\]
Since $k\le j$, we find $m-c\le b-1$, i.e., $d\le m-3$. It follows that 
\[
(A,D)=\brk1{
Z_{1(m-c)}\cup Y_{(m-c)(b-1)}\cup X_{(b-1)(m-2)},\
Y_{1(m-1-c)}\cup Y_{b(m-2)}
}.
\]
Therefore, $N_{41}=\chi(d\le m-3)$.
\item
If $p\in D\times D$, then 
\[
(B,C)=\brk1{X_{0(j-1)},\ Z_{(k+1)(m-1)}}
\quad\text{and}\quad
(j,k)=(b,\,m-1-c).
\]
It follows that 
\[
(A,D)=\brk1{
Z_{1(m-1-c)}\cup Y_{(m-1-c)b}\cup X_{b(m-2)},\
Y_{0(m-2-c)}\cup Y_{(b+1)(m-1)}
}.
\]
Therefore, $N_{42}=1$.

\item\label[itm]{itm:DC}
If $p\in D\times C$, the same argument gives $(j,k)=(b,m-c)$ and $N_{43}=1$.

\item
If $p\in B\times D$, symmetry gives $N_{44}=1$.
\end{enumerate}
Thus $N_4=N_{41}+\dots+N_{44}=3+\chi(d\le m-3)$.
\end{enumerate}
Summing the four cases yields the stated formula for $N=N_1+\dots+N_4$.
\end{proof}

Finally, we count the remaining length-$4$ types containing a part equal to $m$. 

\begin{proposition}\label{prop:mbcd:b<=m-2}
We have the following.
\begin{enumerate}
\item
For $\tau=m(m-2)^22$ with $m\ge 5$, we have $\stc(H_m;\,\tau)=12m^2-8m+54$.
\item
For $\tau=m(m-2)cd$ with $3\le c,d\le m-3$ and $c+d=m$, 
\[
\stc(H_m;\,\tau)
=4c^2
+4dm
+18m+24
+14\binom{m}{d}
+8\binom{m-1}{d-1}
+8\binom{m-2}{d}
+2\sum_{k\in[c, d]}\binom{m-3}{k}.
\]
\item
For $\tau=mbcd$ with $3\le b, c, d\le m-3$ and $b+c+d=2m-2$, 
\begin{multline*}
\stc(H_m;\,\tau)
=24+\sum_{k\in[b,c,d]}\brk4{
4\binom{m}{k}
+10\binom{m-1}{k}
+8\binom{m-2}{k}
+2\binom{m-3}{k}}
\\
+4\brk4{
\binom{b+2}{m-c}
+\binom{c+2}{m-d}
+\binom{d+2}{m-b}
+\binom{b}{m-c-1}
+\binom{c}{m-d-1}
+\binom{d}{m-b-1}
}.
\end{multline*} 
\end{enumerate}
\end{proposition}
\begin{proof} 
For each type $\tau$, we compute $\stc_{23}(\tau)$, $\stc_{24}(\tau)$, and $\stc_{34}(\tau)$ using \cref{lem:s23:mbcd}.
\begin{enumerate}
\item
For $\tau=m(m-2)^22$, we have
\[
\stc_{23}(\tau)
=2m^2-6m+7
\quad\text{and}\quad
\stc_{24}(\tau)
=\stc_{34}(\tau)
=3m+7.
\] 
By \cref{lem:split,lem:s1j:mbcd}, we obtain
\[
\stc(\tau)
=2\brk1{(4m^2-4m+6)
+(2m^2-6m+7)
+2(3m+7)}
=12m^2-8m+54.
\]
\item
For $\tau=m(m-2)cd$,  we have
\begin{align*}
\stc_{23}(\tau)
&=3\binom{m-1}{d}
+\binom{m-3}{d}
+\binom{c+2}{2}
+c+3,
\\
\stc_{24}(\tau)
&=3\binom{m-1}{c}
+\binom{m-3}{c}
+\binom{d+2}{2}
+d+3,
\quad\text{and}\\
\stc_{34}(\tau)
&=2m+4
+\binom{c+2}{2}
+\binom{d+2}{2}.
\end{align*}
By \cref{lem:split,lem:s1j:mbcd}, 
\[
\stc(\tau)
=2\brk4{
m^2+3m-2+4\brk3{\binom{m}{d}+\binom{m-1}{d-1}+\binom{m-2}{d}}
+\stc_{23}(\tau)
+\stc_{24}(\tau)
+\stc_{34}(\tau)},
\]
which simplifies to the desired formula.
\item
For $\tau=mbcd$,  we have
\[
\stc_{23}(\tau)
=4+\binom{m-1}{d}
+2\binom{m-2}{d}
+\binom{m-3}{d}
+\sum_{k\in[b,c]}
\brk4{\binom{k+2}{m-d}
+\binom{k}{m-d-1}}.
\] 
By symmetry, the remaining terms are obtained by summing over the six ordered pairs
\[
(b,c),\ (c,b),\ (b,d),\ (d,b),\ (c,d),\ (d,c).
\]
Consequently,
\begin{align*}
&\stc_{23}(\tau)
+\stc_{24}(\tau)
+\stc_{34}(\tau)
-12
-\sum_{k\in[b,c,d]}
\brk4{
\binom{m-1}{k}
+2\binom{m-2}{k}
+\binom{m-3}{k}}\\
=\ &2\brk4{
\binom{b+2}{m-c}
+\binom{c+2}{m-d}
+\binom{d+2}{m-b}
+\binom{b}{m-c-1}
+\binom{c}{m-d-1}
+\binom{d}{m-b-1}
}.
\end{align*}
Substitution into \cref{lem:split,lem:s1j:mbcd} gives the stated formula.
\end{enumerate}
This completes the proof.
\end{proof}

\begin{theorem}\label{thm:sn-:m3}
For $m\ge 44$, neither $\inc(\mathbf m\times \mathbf 3)$ nor $H_m$ is strongly nice. Moreover,
\begin{align*}
\stc\brk1{\inc(\mathbf m\times\mathbf 3);\, \lambda}
&=\stc\brk1{H_m;\,\lambda^-}
=m^4-\frac{25}{3}m^3+36m^2-\frac{194}{3}m+144,\quad\text{and}\\
\stc\brk1{\inc(\mathbf m\times\mathbf 3);\, \mu}
&=\stc\brk1{H_m; \mu^-}
=m^4-\frac{26}{3}m^3+52m^2-\frac{379}{3}m+240,
\end{align*}
where 
\begin{align*}
&\lambda=(m+2)(m-2)(m-4)4, \quad
\lambda^-=m(m-2)(m-4)4, \\
&\mu=(m+2)(m-3)^24, 
\quad\text{and}\quad
\mu^-=m(m-3)^24.
\end{align*}
\end{theorem}
\begin{proof}
Let $G=\inc(\mathbf m\times\mathbf 3)$. Recall that $G$ is the disjoint union of $H_m$ and two isolated vertices. Since~$\alpha(H_m)=m$, every stable block of size $m+2$ in $G$ consists of a stable block of size $m$ in $H_m$ together with the two isolated vertices. Therefore,
\[
\stc(G;\lambda)=\stc(H_m;\lambda^-)
\qquad\text{and}\qquad
\stc(G;\mu)=\stc(H_m;\mu^-).
\]
Applying Part~(2) of \cref{prop:mbcd:b<=m-2} with $(c,\,d)=(m-4,\,4)$ to
$\lambda^-$
gives the first formula in the statement. Applying Part~(3) of \cref{prop:mbcd:b<=m-2} with
\[
(b,\,c,\,d)=(m-3,\,m-3,\,4)
\]
to $\mu^-$
gives the second formula in the statement. Set
\[
D(m)=\stc(H_m;\lambda^-)-\stc(H_m;\mu^-).
\]
Subtracting the two formulas gives
\[
D(m)=\frac{1}{3}\bigl(m^3-48m^2+185m-288\bigr).
\]
Moreover, $D(44)=36>0$, and for $m\ge 44$, 
\[
D(m+1)-D(m)=m^2-31m+46=m(m-31)+46>0.
\]
Hence $D(m)>0$ for every $m\ge44$. Since $\lambda$ strictly dominates $\mu$ and $\lambda^-$ strictly dominates $\mu^-$, neither $G$ nor $H_m$ is strongly nice.
\end{proof}

In view of \cref{prop:s+:sn+:n+}, \cref{thm:sn-:m3} implies that the graph $\inc(\mathbf m\times\mathbf 3)$ is not Schur positive for~$m\ge 44$. We now strengthen this conclusion to every $m\ge7$.

\begin{theorem}\label{thm:ns:m3}
For every positive integer $m$, the lattice $\mathbf m\times\mathbf 3$ is not Schur positive if and only if~$m\ge7$. More precisely, for $m=7$, one has $[s_{(6,5,5,5)}]X_{\inc(\mathbf 7\times\mathbf 3)}=-36$, and for every $m\ge 8$,
\[
[s_{(m+2)(m-3)^24}]X_{\inc(\mathbf m\times\mathbf 3)}
=-\frac{4}{3}m^3+16m^2-\frac{176}{3}m+96<0.
\]

\end{theorem}
\begin{proof}
Fix $m\ge8$ and let $\mu=(m+2)(m-3)^24$. Since $s_1^2=s_2+s_{11}$, \cref{prop:Pieri} shows that the partitions $\eta$ satisfying $[s_\mu]s_\eta s_1^2\ne0$ are obtained from $\mu$ by removing a horizontal or vertical strip of size $2$. Among these shapes, only $\mu^-=m(m-3)^24$ contributes to $[s_\mu]X_{\inc(\mathbf m\times\mathbf3)}$. Indeed, every other such partition $\eta$ satisfies $\eta_1\ge m+1$. In any special ribbon tabloid of shape $\eta$, the ribbon containing the rightmost cell has its head in the first column and is connected. It therefore contains at least one cell in every intervening column and has size at least $\eta_1\ge m+1$. Since $\alpha(H_m)=m$, we have $\stc(H_m;\lambda_T)=0$ for every $T\in\mathcal T_\eta$.
\begin{figure}[h]
\begin{tikzpicture}[scale=0.45, round]
\input H3-sh.tex
\draw[ribbon] (.5,  .5)--(1.5,  .5);
\draw[ribbon] (.5, 1.5)--(2.5, 1.5)--(2.5,  .5)--(5.5,  .5);
\draw[ribbon] (.5, 2.5)--(3.5, 2.5)--(3.5, 1.5)--(6.5, 1.5)--(6.5, .5);
\foreach \x/\y in {5.5/.5, 1.5/.5}{\fill[black] (\x,\y) circle (3pt);}
\node at (kappa) {$\kappa_1=(m-1)(m-3)2m$};

\begin{scope}[xshift=11cm]
\input H3-sh.tex
\draw[ribbon] (.5,  .5)--(5.5,  .5);
\draw[ribbon] (.5, 1.5)--(2.5, 1.5);
\draw[ribbon] (.5, 2.5)--(3.5, 2.5)--(3.5, 1.5)--(6.5, 1.5)--(6.5, .5);
\foreach \x/\y in {5.5/.5, 2.5/1.5}{\fill[black] (\x,\y) circle (3pt);}
\node at (kappa) {$\kappa_2=(m-1)3(m-4)m$};
\end{scope}

\begin{scope}[xshift=22cm]
\input H3-sh.tex
\draw[ribbon] (.5,  .5)--(1.5,  .5);
\draw[ribbon] (.5, 1.5)--(2.5, 1.5)--(2.5,  .5)--(6.5,  .5);
\draw[ribbon] (.5, 2.5)--(3.5, 2.5)--(3.5, 1.5)--(6.5, 1.5);
\foreach \x/\y in {1.5/.5, 6.5/1.5}{\fill[black] (\x,\y) circle (3pt);}
\node at (kappa) {$\kappa_3=(m-2)^22m$};
\end{scope}

\begin{scope}[yshift=-6cm]
\input H3-sh.tex
\draw[ribbon] (.5,  .5)--(6.5,  .5);
\draw[ribbon] (.5, 1.5)--(2.5, 1.5);
\draw[ribbon] (.5, 2.5)--(3.5, 2.5)--(3.5, 1.5)--(6.5, 1.5);
\foreach \x/\y in {2.5/1.5, 6.5/1.5}{\fill[black] (\x,\y) circle (3pt);}
\node at (kappa) {$\kappa_4=(m-2)3(m-3)m$};
\end{scope}

\begin{scope}[yshift=-6cm, xshift=11cm]
\input H3-sh.tex
\draw[ribbon] (.5,  .5)--(5.5,  .5);
\draw[ribbon] (.5, 1.5)--(6.5, 1.5)--(6.5, .5);
\draw[ribbon] (.5, 2.5)--(3.5, 2.5);
\foreach \x/\y in {3.5/2.5, 5.5/.5}{\fill[black] (\x,\y) circle (3pt);}
\node at (kappa) {$\kappa_5=4(m-2)(m-4)m$};
\end{scope}

\begin{scope}[yshift=-6cm, xshift=22cm]
\input H3-sh.tex
\draw[ribbon] (.5,  .5)--(6.5,  .5);
\draw[ribbon] (.5, 1.5)--(6.5, 1.5);
\draw[ribbon] (.5, 2.5)--(3.5, 2.5);
\foreach \x/\y in {3.5/2.5, 6.5/1.5}{\fill[black] (\x,\y) circle (3pt);}
\node at (kappa) {$\kappa_6=4(m-3)^2m$};
\end{scope}
\end{tikzpicture}
\caption{The six contributing special ribbon tabloids
$\kappa_1,\ldots,\kappa_6$ in $\mathcal T_{m(m-3)^2 4}$.}
\label{fig:m.m-3.m-3.4}
\end{figure}
A special ribbon tabloid in $\mathcal T_{\mu^-}$ can contribute only if every part of its ribbon content is at most $\alpha(H_m)=m$. Hence it suffices to consider the six tabloids $\kappa_1,\dots,\kappa_6$ shown in \cref{fig:m.m-3.m-3.4}. Among these, the signs of $\kappa_2$, $\kappa_3$, and
$\kappa_6$ are positive, whereas the signs of $\kappa_1$, $\kappa_4$, and
$\kappa_5$ are negative. By \cref{prop:m.m-1.c.d,prop:mbcd:b<=m-2}, 
\begin{align*}
&[s_\mu]X_{\inc(\mathbf m\times\mathbf 3)}
=[s_{\mu^-}]X_{H_m}
\\
=&-\stc(H_m; \kappa_1)
+\stc(H_m; \kappa_2)
+\stc(H_m; \kappa_3)
-\stc(H_m; \kappa_4)
-\stc(H_m; \kappa_5)
+\stc(H_m; \kappa_6)\\
=&-(9m^2-19m+54)
+(3m^3-18m^2+41m+18)
\\
&+(12m^2-8m+54)
-(4m^3-15m^2+49m+18)\\
&-\brk2{m^4-\frac{25}{3}m^3+36m^2-\frac{194}{3}m+144} 
+\brk2{m^4-\frac{26}{3}m^3+52m^2-\frac{379}{3}m+240}\\
=&-\frac{4}{3}m^3+16m^2-\frac{176}{3}m+96\\
=&-\frac{4}{3}
\brk2{(m-8)^3+12(m-8)^2+44(m-8)+24}<0.
\end{align*}

For $m=7$, the identity $X_{\inc(\mathbf7\times\mathbf3)}=s_1^2X_{H_7}$, \cref{prop:Pieri}, and the exact computation in \cref{app:small-cases} give
\[
\begin{aligned}
[s_{(6,5,5,5)}]X_{\inc(\mathbf 7\times\mathbf 3)}
&=
[s_{(6,5,5,3)}]X_{H_7}
+[s_{(6,5,4,4)}]X_{H_7}
+2[s_{(5,5,5,4)}]X_{H_7}  \\
&=108-12+2(-66)=-36.
\end{aligned}
\]
Thus $\mathbf7\times\mathbf3$ is not Schur positive.

The exact computations described in \cref{app:small-cases} show that $\mathbf m\times\mathbf3$ is Schur positive for $m\le6$. This completes the proof.

\end{proof}

\section*{Acknowledgments}
The authors thank the anonymous referees for their careful reading and helpful comments, which improved the exposition.

\appendix

\section{Exact verification of the finite cases}
\label{app:small-cases}
All computations in this appendix were performed in SageMath~$10.9$ using exact integer arithmetic. For each pair $(m,n)$ listed in \cref{tab:finite-verification}, let
\[
G=\inc(\mathbf m\times\mathbf n).
\]
\begin{table}[ht]
\centering
\caption{Exact verification of the finite cases.}
\label{tab:finite-verification}
\small
\begin{tabular}{@{}cccc@{}}
\toprule
Lattice
& $m$
& Number of partitions $\mu\vdash mn$
& Number of nonzero Schur coefficients\\
\midrule
$\mathbf m\times\mathbf 2$
& 1--7
& $2,\,5,\,11,\,22,\,42,\,77,\,135$
& $2,\,4,\,9,\,18,\,35,\,65,\,115$\\
$\mathbf m\times\mathbf 3$
& 1--6
& $3,\,11,\,30,\,77,\,176,\,385$
& $3,\,9,\,22,\,56,\,128,\,283$\\
\bottomrule
\end{tabular}
\end{table}
The Schur coefficients of $X_G$ were computed directly from
\cref{thm:scoeff}, rather than by converting from another basis.
For each partition $\mu\vdash mn$, the program performs the following steps:
\begin{enumerate}
\item it generates the set $\mathcal T_\mu$ of all special ribbon
tabloids of shape $\mu$;

\item for each $T\in\mathcal T_\mu$, it determines the sign
$\sgn(T)$ and the partition $\lambda_T$ underlying the ribbon
content of $T$;

\item it enumerates the stable compositions of $G$ of type
$\lambda_T$ and computes
\[
\stc(G;\lambda_T)
=\abs{\mathcal{SC}(G;\lambda_T)};
\]

\item it evaluates the signed sum
\[
[s_\mu]X_G
=
\sum_{T\in\mathcal T_\mu}
\sgn(T)\stc(G;\lambda_T).
\]
\end{enumerate}
Repeating this procedure for every partition $\mu\vdash mn$
produces the complete Schur expansion of $X_G$. In particular,
Schur positivity in the finite cases is verified by checking that
every coefficient obtained from the displayed signed sum is
nonnegative. 

For the boundary case $(m,n)=(7,3)$, the same procedure was applied
to the reduced incomparability graph $H_7$. It gives
\[
[s_{(6,5,5,3)}]X_{H_7}=108,
\quad
[s_{(6,5,4,4)}]X_{H_7}=-12,
\quad\text{and}\quad
[s_{(5,5,5,4)}]X_{H_7}=-66.
\]
Consequently, by \cref{prop:Pieri},
\[
\begin{aligned}
[s_{(6,5,5,5)}]
X_{\inc(\mathbf7\times\mathbf3)}
&=
[s_{(6,5,5,3)}]X_{H_7}
+[s_{(6,5,4,4)}]X_{H_7}
+2[s_{(5,5,5,4)}]X_{H_7}\\
&=108-12+2(-66)=-36.
\end{aligned}
\]
\Cref{tab:finite-verification} summarizes the complete finite verification. Within each row, the entries are listed in increasing order of $m$.

\bibliographystyle{abbrvnat}
		
\bibliography{csf.bib}
\end{document}

%% file: H2-sh.tex
\draw[gray] (0,1)--(7,1);
\draw[gray] (0,2)--(4,2);
\foreach \x in {1, 2, 3, 4, 5, 6} {\draw[gray] (\x, 0)--(\x, 2);}
\foreach \x in {1, 2, 3} {\draw[gray] (\x, 2)--(\x, 3);}
\draw[border] (0,0)--(7,0)--(7,2)--(4,2)--(4,3)--(0,3)--cycle;
\node at (4.56, 1.7) {\ldots};
\node at (4.56, 0.7) {\ldots};
\foreach \x/\y in {0.5/2.5, 0.5/1.5, 0.5/0.5, 6.5/0.5}{
  \fill[black] (\x, \y) circle (3pt);}
\coordinate (kappa) at (3.5, -.8);

%% file: H2-rm.tex
\draw[gray] (0,1)--(6,1);
\draw[gray] (0,2)--(4,2);
\foreach \x in {1, 2, 3, 4, 5} {
  \draw[gray] (\x, 0)--(\x, 2);}
\draw[gray] (1, 2)--(1, 3);
\draw[gray] (2, 2)--(2, 3);
\node at (4.56, 1.5) {\ldots};
\node at (4.56, 0.5) {\ldots};

%% file: H3-sh.tex
\draw[gray] (0,0) -- (7,0);
\draw[gray] (0,1) -- (7,1);
\draw[gray] (0,2) -- (4,2);
\foreach \x in {1,2,3,4,5,6,7} {\draw[gray] (\x, -1)--(\x, 2);}
\foreach \x in {8,9} {\draw[gray] (\x, -1)--(\x, 0);}
\foreach \x in {1,2,3} {\draw[gray] (\x, 2)--(\x, 3);}
\foreach \x in {1.7, .7, -.3} {\node at (4.56, \x) {\dots};}
\foreach \x/\y in {.5/-.5, .5/.5, .5/1.5, .5/2.5, 6.5/.5, 9.5/-.5}{\fill[black] (\x,\y) circle (3pt);}
\draw[ribbon] (.5, -.5)--(9.5, -.5);
\draw[border] (0, -1)--(10, -1)--(10, 0)--(7,0)--(7,2)--(4,2)--(4,3)--(0,3)--cycle;
\coordinate (kappa) at (5, -1.8);

%% file: csf.bib
@misc{Sv26X,
	author = {Ethan Shelburne and Stephanie van Willigenburg},
	howpublished = {arXiv:2604.26158},
	title = {A {Schur-positivity} classification for complete multipartite graphs},
	year = {2026}}

@article{LQYZ25,
	author = {Li, Grace M. X. and Qiu, Dun and Yang, Arthur L. B. and Zhang, Zhong-Xue},
	journal = {SIAM Discrete Math.},
	number = {4},
	pages = {2250--2267},
	title = {Stanley's conjecture on the {Schur} positivity of distributive lattices},
	volume = {39},
	year = {2025}}

@article{Ale21,
	author = {Alexandersson, Per},
	journal = {J. Alg. Combin.},
	number = {2},
	pages = {299--325},
	title = {{LLT polynomials}, elementary symmetric functions and melting lollipops},
	volume = {53},
	year = {2021}}

@article{LE99,
	author = {Zbigniew Lonc and Muktar E. Elzobi},
	journal = {J. Combin. Theory Ser. A},
	number = {1},
	pages = {140--150},
	title = {Chain Partitions of Products of Two Chains},
	volume = {86},
	year = {1999}}

@article{Gri98,
	author = {Jerrold R. Griggs},
	journal = {Discrete Math.},
	number = {1},
	pages = {157-162},
	title = {Problems on chain partitions},
	volume = {72},
	year = {1988}}

@book{BirGarrett1940B,
	address = {Providence, RI},
	author = {Garrett Birkhoff},
	edition = {3rd},
	publisher = {Amer. Math. Soc.},
	series = {American Mathematical Society colloquium publications},
	title = {Lattice Theory},
	volume = {25},
	year = {1967}}

@article{TV26-vglue,
	author = {Foster Tom and Aarush Vailaya},
	journal = {Electron. J. Combin.},
	number = {2},
	pages = {\#P2.9},
	title = {The chromatic symmetric function of graphs glued at a single vertex},
	volume = {33},
	year = {2026}}

@article{WZ25,
	author = {Wang, David G. L. and Zhou, James Z.F.},
	journal = {Adv. Appl. Math.},
	pages = {102886},
	publisher = {Elsevier BV},
	title = {A composition method for neat formulas of chromatic symmetric functions},
	volume = {167},
	year = {2025}}

@inproceedings{Erd84,
	address = {Berlin},
	author = {P. Erd{\H o}s},
	booktitle = {Measure Theory: Proceedings of the Conference held at Oberwolfach, June 26--July 2, 1983},
	editor = {K{\"o}lzow, Dietrich and Maharam-Stone, Dorothy},
	pages = {321--327},
	publisher = {Springer-Verlag},
	series = {Lecture Notes in Math. (LNM) 1089},
	title = {Some combinatorial, geometric and set theoretic problems in measure theory},
	year = {1984}}

@article{Sv25,
	author = {Ethan Shelburne and Stephanie van Willigenburg},
	journal = {Enumer. Combin. Appl.},
	number = {1},
	pages = {Article S2R8},
	title = {Schur-positivity for generalized nets},
	volume = {5},
	year = {2025}}

@article{Hik25,
	author = {Tatsuyuki Hikita},
	journal = {S\'{e}m. Lothar. Combin.},
	number = {Article \#31},
	pages = {12 pp.},
	title = {On the {Stanley--Stembridge} conjecture},
	volume = {93B},
	year = {2025}}

@article{LLYZ25,
	author = {Ethan Y. H. Li and Grace M. X. Li and Arthur L. B. Yang and Zhong-Xue Zhang},
	journal = {S\'{e}m. Lothar. Combin.},
	pages = {Art. \#154, 11 pp.},
	title = {Strongly nice property and {Schur} positivity of graphs},
	volume = {93B},
	year = {2025}}

@article{Bir1912,
	author = {George D. Birkhoff},
	journal = {Ann. of Math.},
	number = {1/4},
	pages = {42--46},
	publisher = {Annals of Mathematics},
	title = {A Determinant Formula for the Number of Ways of Coloring a Map},
	volume = {14},
	year = {1912}}

@article{ER90,
	author = {{\"O}mer E{\u g}ecio{\u g}lu and Jeffrey B. Remmel},
	journal = {Linear Multilinear Algebra},
	number = {1-2},
	pages = {59--84},
	title = {A combinatorial interpretation of the inverse {K}ostka matrix},
	volume = {26},
	year = 1990}

@misc{TW23X,
	author = {Thibon, Jean-Yves and Wang, David G. L.},
	howpublished = {arXiv:2305.07858},
	title = {A noncommutative approach to the {S}chur positivity of chromatic symmetric functions},
	year = {2023}}

@article{Sta81,
	author = {Richard P. Stanley},
	journal = {Bull. Amer. Math. Soc.},
	number = {2},
	pages = {254--265},
	title = {{Review: I. G. Macdonald, Symmetric functions and Hall polynomials}},
	volume = {4},
	year = {1981}}

@article{WW20,
	author = {Wang, David G. L. and Wang, Monica M. Y.},
	journal = {Discrete Appl. Math.},
	pages = {621--630},
	title = {A combinatorial formula for the {S}chur coefficients of chromatic symmetric functions},
	volume = {285},
	year = {2020}}

@article{GS01,
	author = {Gebhard, David D. and Sagan, Bruce E.},
	journal = {J. Alg. Combin.},
	number = {3},
	pages = {227--255},
	title = {A chromatic symmetric function in noncommuting variables},
	volume = {13},
	year = {2001}}

@book{Sag01B,
	author = {Sagan, Bruce E.},
	edition = {2nd},
	publisher = {Springer-Verlag, New York},
	series = {Grad. Texts in Math.},
	title = {The Symmetric Group: Representations, Combinatorial Algorithms, and Symmetric Functions},
	volume = {203},
	year = {2001}}

@book{Mac95B,
	author = {I. G. Macdonald},
	edition = {2nd},
	publisher = {The Clarendon Press, Oxford Univ. Press, New York},
	series = {Oxford Math. Monogr.},
	title = {Symmetric {F}unctions and {H}all {P}olynomials},
	year = {1995}}

@article{SS93,
	author = {Stanley, Richard P. and Stembridge, John R.},
	bdsk-color = {4},
	journal = {J. Combin. Theory Ser. A},
	number = {2},
	pages = {261--279},
	title = {On immanants of {Jacobi--Trudi} matrices and permutations with restricted position},
	volume = {62},
	year = {1993}}

@article{Sta95,
	author = {Stanley, Richard P.},
	bdsk-color = {4},
	journal = {Adv. Math.},
	number = {1},
	pages = {166--194},
	title = {A symmetric function generalization of the chromatic polynomial of a graph},
	volume = {111},
	year = {1995}}

@article{Dv18,
	author = {Dahlberg, Samantha and {van Willigenburg}, Stephanie},
	journal = {SIAM J. Discrete Math.},
	number = {2},
	pages = {1029--1039},
	title = {Lollipop and lariat symmetric functions},
	volume = {32},
	year = {2018}}

@article{Gas99,
	author = {Gasharov, Vesselin},
	journal = {Discrete Math.},
	number = {1--3},
	pages = {229--234},
	title = {On {S}tanley's chromatic symmetric function and clawfree graphs},
	volume = {205},
	year = {1999}}

@article{Sta98,
	author = {Stanley, Richard P.},
	journal = {Discrete Math.},
	number = {1-3},
	pages = {267--286},
	title = {Graph colorings and related symmetric functions: ideas and applications: a description of results, interesting applications, \& notable open problems},
	volume = {193},
	year = {1998}}

@article{Gas96,
	author = {Gasharov, Vesselin},
	bdsk-color = {4},
	booktitle = {Proceedings of the 6th {C}onference on {F}ormal {P}ower {S}eries and {A}lgebraic {C}ombinatorics ({N}ew {B}runswick, {NJ}, 1994)},
	journal = {Discrete Math.},
	number = {1--3},
	pages = {193--197},
	title = {Incomparability graphs of {$(3+1)$}-free posets are {$s$}-positive},
	volume = {157},
	year = {1996}}

@book{Sta99B,
	author = {Stanley, Richard P.},
	pages = {xii+581},
	publisher = {Camb. Univ. Press, Cambridge},
	series = {Cambridge Stud. in Adv. Math.},
	title = {Enumerative Combinatorics. {V}ol. 2},
	volume = {62},
	year = {1999}}

@article{SW16,
	author = {Shareshian, John and Wachs, Michelle L.},
	journal = {Adv. Math.},
	pages = {497--551},
	title = {Chromatic quasisymmetric functions},
	volume = {295},
	year = {2016}}
